\documentclass[11pt]{amsart}
\allowdisplaybreaks[1]

\usepackage{a4wide}
\usepackage{amssymb}
\usepackage{amsmath}
\usepackage{amscd} 

\newtheorem{theorem}{Theorem}
\newtheorem{lemma}{Lemma}
\newtheorem{prop}[lemma]{Proposition}
\newtheorem{coro}[lemma]{Corollary}

\newtheorem{definition}[lemma]{Definition}

\numberwithin{lemma}{section}
\numberwithin{theorem}{section}
\numberwithin{fact}{section}

\newcommand{\ts}{\hspace{0.5pt}}

\newcommand{\CC}{\mathbb{C}\ts}
\newcommand{\RR}{\mathbb{R}\ts}

\newcommand{\NN}{\mathbb{N}}
\newcommand{\TT}{\mathbb{T}}

\newcommand{\oplam}{\mbox{\Large $\curlywedge$}}

\newcommand{\dd}{\,{\rm d}}

\newcommand{\MM}{\mathcal{M}(G)}
\newcommand{\MCV}{\mathcal{M}_{C,V}(G)}
\newcommand{\MTB}{\mathcal{M}^{\infty}(G)}

\newcommand{\cM}{\mathcal{M}}

\newcommand{\Oomega}{(\varOmega,\alpha)}

\newcommand{\Ttheta}{(\varTheta,\beta)}

\newcommand{\LO}{L^2 (\varOmega,m)}

\newcommand{\cT}{\mathcal{T}}

\newcommand{\Ghat}{\widehat{G}}
\newcommand{\gammahat}{\widehat{\gamma}}

\newcommand{\supp}{\mbox{supp}}

\newcommand{\Hm}[1]{\leavevmode{\marginpar{\tiny%
$\hbox to 0mm{\hspace*{-0.5mm}$\leftarrow$\hss}%
\vcenter{\vrule depth 0.1mm height 0.1mm width \the\marginparwidth}%
\hbox to
0mm{\hss$\rightarrow$\hspace*{-0.5mm}}$\\\relax\raggedright #1}}}

\newcommand{\cp}{(G,H,\widetilde{L})}
\newcommand{\mcp}{(G,H,\widetilde{L},\rho)}

\begin{document}

\title[Pure point diffraction  ]{ Pure point diffraction and cut and project schemes for measures:  The smooth case}

\author{Daniel Lenz}
\address{Fakult\"at f\"ur Mathematik, TU Chemnitz,
09107 Chemnitz, Germany}
\email{dlenz@mathematik.tu-chemnitz.de}

\author{Christoph Richard}
\address{Fakult\"at f\"{u}r Mathematik, Universit\"{a}t Bielefeld,
Postfach 100131, 33501 Bielefeld, Germany}
\email{richard@math.uni-bielefeld.de}

\begin{abstract} 
We present cut and project formalism based on measures and
  continuous weight functions of sufficiently fast decay. The emerging
  measures are strongly almost periodic. The corresponding dynamical
  systems are compact groups and homomorphic images of the underlying
  torus. In particular, they are strictly ergodic with pure point
  spectrum and continuous eigenfunctions. Their diffraction can be
  calculated explicitly. Our results cover and extend corresponding 
  earlier results on dense Dirac combs and continuous weight functions 
  with compact support. They also mark a clear difference in terms of 
  factor maps between the case of continuous and non-continuous weight
  functions. 
\end{abstract}

\maketitle

\section{Introduction}

This paper is concerned with the harmonic analysis behind certain
models of aperiodic order.  The latter is a specific form of order
with long range correlations but no translation symmetry.  It has
attracted a lot of attention in the last two decades, compare the surveys and
monographs \cite{BMed,Jan,Moo,Pat,Sen}. 

This attention is partly due to the actual discovery of physical
substances, later called quasicrystals, exhibiting such a form of
order \cite{SBGC,INF}. Their key feature is a pure point diffraction
spectrum combined with a non-periodic structure. (In a periodic
structure pure point diffraction results easily from a Poisson
summation type formula, see \cite{Cord,Lag} for further ideas in this
direction.) This attention is also due to the conceptual mathematical
relevance of aperiodic order as an intermediate form of (dis)order
between periodicity and randomness. In fact, aperiodic order 
has highly distinctive and far from being understood geometric,
combinatorial and Fourier analytic features.

The most prominent models of aperiodic order arise from so called cut
and project schemes. They are called model sets or harmonious
sets. They were introduced and first studied by Meyer in \cite{Mey}
for purely theoretical reasons. His investigations have later been
generalised and extended in various directions (see
\cite{Moody-old,Moody} for recent surveys and \cite{Rob} for a recent 
inverse spectral type result).  In the physics community
cut and project models have been the objects of choice from the very
beginning of theoretical investigation of quasicrystals
\cite{LS}.

\smallskip

In the study of aperiodic order and diffraction the use of dynamical systems
has a long history going back to \cite{Dwo,Que} (see
\cite{Hof,Martin2,Boris,Boris2} as well).  Recently two further lines
of research have proven fruitful: These are the systematic studies of
notions of almost periodicity \cite{BM,Gouere-2,MS} and the
replacement of sets by translation bounded measures
\cite{BL,BL2,Ric,LMS-1}.

In line with these developments the basic aim of this paper is to
extend the cut and project formalism to measures. More precisely, specific goals  of this paper are

\begin{itemize}

\item to develop a cut and project scheme based on measures (instead of sets),

\item to study the dynamical systems arising from these schemes,

\item to investigate almost periodicity properties in this context.

\end{itemize}

Our results lead to a rather complete picture with quite strong
properties being valid, provided the weight function is
\textit{smooth}, i.e., continuous and sufficiently fast decaying. In
this case, almost periodicity is present in a rather strong form and
(essentially) everything is determined by the underling torus
dynamical system. More precisely, the arising measure dynamical
systems are factors of the torus dynamical system.  They carry a group
structure and the factor map is a group homomorphism. Diffraction can
be calculated explicitly.

In some sense our  models are more regular than the ``usual'' cut and
project schemes, where the weight function is the characteristic
function of a Riemann integrable set. Mathematically, this is
reflected in the almost periodicity properties of the
underlying measures (as opposed to almost periodicity properties of
averaged quantities like the autocorrelation). From the point of view
of physics one may also argue in favour of our models: The strict cut
off procedure in the usual model sets is highly idealised, whereas the
cut off by continuous functions may be more realistic, at least in an
averaged sense. Moreover, such models are used to analyse diffraction 
properties of random tilings \cite{E,Hen}, whose vertex sets are 
derived from model sets, see also \cite{BMRS,Ric}.

\smallskip

For special cases some of these results are already known.  Hof \cite{Hof}
presents results on continuous weight functions with compact support as a tool
in his study of the usual model sets. Richard \cite{Ric} has systematically
investigated dense Dirac combs on $\RR^d$.  Our results cover and considerably
extend the corresponding results of these authors, see Section \ref{Dense}
below.

We would like to emphasise that these results do no longer hold if
the smoothness of the weight function is violated. More precisely, in
the usual cut and project schemes,  the arising
dynamical system is neither a group nor a factor of the torus. On the
contrary, the torus in that case  is a factor of the dynamical system
\cite{BHP,Martin2}, but not vice versa (see below Section
\ref{Complementary}).  Thus, our results show in particular a change in the
role of the torus system depending on the continuity of the weight
function. 

\smallskip

The paper is organised as follows: In Section \ref{Measure} we recall
background and notation. Section \ref{Cut} presents the cut and
project schemes for measures and gives our main results. The necessary
investigation of almost-periodicity is carried out in Section
\ref{Strongly}. An abstract study of factors in our context is given
in Section \ref{Factors}.  Section \ref{Study} studies the dynamical
systems arising from the measure cut and project schemes.  After these
preparations we discuss the proof of Theorem \ref{main} in Section 
\ref{Proof_main}.  A Weyl formula on uniform
distribution is presented in Section \ref{Weyl}.  This is used to
discuss the so-called Fourier Bohr coefficients and the proof of
Theorem \ref{main_diffraction} in Section \ref{Fourier}.  Dense Dirac
combs and other examples are studied in Section \ref{Dense}.
Injectivity of the arising factor map is discussed in Section
\ref{Injectivity}. Finally, in Section \ref{Complementary}, we compare
our results to those for usual model set dynamical systems.

\section{Measure dynamical systems and diffraction}\label{Measure} 

\subsection{Dynamical systems} Whenever $X$ is a $\sigma$-compact locally
compact space  
(by which we mean to include
the Hausdorff property),  the space of continuous functions on
$X$ is denoted by $C(X)$, the subspace of continuous functions with compact
support by $C_c (X)$ and the space of continuous bounded functions by $C_b
(X)$. The latter two spaces are complete normed spaces when equipped with the 
supremum norm $\|\cdot\|_\infty$.  

A topological space $X$ carries the Borel $\sigma$-algebra generated
by all closed subsets of $X$.  By the Riesz-Markov representation
theorem, the set $\mathcal{M} (X)$ of all complex regular Borel
measures on $X$ can then be identified with the dual space $C_c
(X)^\ast$ of complex valued, linear functionals on $C_c(X)$ which are
continuous with respect to a suitable topology, see \cite[Ch.\
6.5]{Ped} for details. In particular, we write $\int_X f \dd\mu =
\mu(f)$ for $f\in C_c(X)$. The space $\mathcal{M} (X)$ carries the
vague topology, i.e., the weakest topology that makes all functionals
$\mu\mapsto \mu(\varphi)$, $\varphi\in C_c (X)$, continuous.
Alternatively, the vague topology arises by considering $\mathcal{M}
(X)$ to be a subset of $\prod_{\varphi \in C_c (X)} \CC$, which is
equipped with the product topology, via
$$ \mathcal{M} (X)\longrightarrow \prod_{\varphi \in C_c (X)} \CC, \;\: \nu
\mapsto (\varphi \mapsto \nu(\varphi)).$$
The total variation of a measure $\mu \in \mathcal{M} (X)$ 
is denoted by $|\mu|$.

\smallskip

We will have to deal with various abelian groups. The group operation
will be written additively as $+$ or $\dotplus$ if necessary to avoid
misunderstandings. Now, let $G$ be a $\sigma$-compact locally compact
abelian (LCA) group. The Haar measure on $G$ is denoted by $m_G$ or
${\rm d} t$. The dual group of $G$ is denoted by $\Ghat$, and the
pairing between a character $\hat{s} \in \Ghat$ and $t\in G$ is
written as $(\hat{s},t)$. As usual the Fourier transform $\widehat{f}$
of an integrable function $f$ is defined by $\widehat{f} (\hat{s}) =
\int_G \overline{(\hat{s},t)} f(t) {\rm d}t$.

Whenever $G$ 
acts on the compact space $\varOmega$ (which is then also Hausdorff by our
convention) by a continuous action
\begin{equation*}
   \alpha \! : \; G\times \varOmega \; \longrightarrow \; \varOmega
   \, , \quad (t,\omega) \, \mapsto \, \alpha^{}_{t} (\omega) \, ,
\end{equation*} 
where $G\times \varOmega$ carries the product topology, the pair 
$\Oomega$ is called a {\em topological dynamical system\/} over $G$. 

An $\alpha$-invariant probability measure on $\varOmega$ is then
called {\em ergodic\/} if every measurable invariant subset of
$\varOmega$ has either measure zero or measure one. The dynamical
system $\Oomega$ is called {\em uniquely ergodic\/} if there exists a
unique $\alpha$-invariant probability measure on $\varOmega$, which
then is ergodic by standard theory. $\Oomega$ is called {\em
minimal\/} if, for all $\omega\in\varOmega$, the $G$-orbit
$\{\alpha^{}_t\ts \omega : t \in G\}$ is dense in $\varOmega$. If
$\Oomega$ is both uniquely ergodic and minimal, it is called {\em
strictly ergodic}.

\smallskip

Given an $\alpha$-invariant probability measure $m$ on $\varOmega$ ,
we can form the Hilbert space $\LO$ of square integrable measurable
functions on $\varOmega$. This space is equipped with the inner
product
\begin{equation*}
    \langle f, g\rangle \; = \;
    \langle f, g\rangle^{}_\varOmega \; := \; 
    \int_\varOmega \overline{f(\omega)}\, g(\omega) \dd m(\omega).
\end{equation*}
The action $\alpha$ gives rise to a unitary representation $T :=
T^\varOmega := T^{(\varOmega,\alpha,m)}$ of $G$ on $\LO$ by
\begin{equation*}
  T_t \! : \; \LO \; \longrightarrow \; \LO \, , 
  \quad (T_t f) (\omega) \; := \;
  f(\alpha^{}_{-t}\ts \omega) \, ,
\end{equation*}
for every $f\in \LO$ and arbitrary $t\in G$. 
 An $f\in \LO$ is called an {\em eigenfunction\/} of $T$ with
{\em eigenvalue\/} $\hat{s}\in \Ghat$ if $T_t f = (\hat{s}, t) f$ for every
$t\in G$.  An eigenfunction (to $\hat{s}$, say) is called {\em continuous\/} 
if it has a continuous representative $f$  with
$f(\alpha^{}_{-t} \ts \omega) = (\hat{s},t)\ts f(\omega)$, for all
$\omega\in\varOmega$ and $t\in G$. The representation $T$ is said to have 
{\em pure point spectrum\/} if the set of eigenfunctions is total in $\LO$.  
One then also says that the dynamical system $\Oomega$ has
{\em pure point dynamical spectrum}.

Finally, we will need the notion of factor of a dynamical system. 

\begin{definition}  
  Let two topological dynamical systems\/ $\Oomega$ and\/ $\Ttheta$
  under the action of $G$ be given.  Then, $\Ttheta$ is called a\/
  {\em factor} of $\Oomega$, with factor map\/ $\varPhi$, if\/
  $\varPhi \! : \varOmega \longrightarrow\varTheta $ is a continuous
  surjection with $\varPhi (\alpha^{}_t (\omega)) = \beta^{}_t
  (\varPhi (\omega))$ for all\/ $\omega\in \varOmega$ and\/ $t\in G$.
\end{definition}

\subsection{Measure dynamical systems}

We will be concerned with dynamical systems built from measures.  These
systems will be discussed next. They have been introduced in \cite{BL,BL2}, to
which we refer for further details and proofs of the subsequent discussion.

\smallskip

A measure $\nu\in \MM $ is called {\em translation bounded\/} if
there exist some  $C>0$ and an open non empty relatively compact  set $V$ in $G$  so
that
\begin{equation}\label{mcv} |\nu| (t+V) \leq C
\end{equation}
for every $t\in G$, where $|\nu|$ is the total variation measure of $\nu$.  
The set of all translation bounded measures
satisfying \eqref{mcv} is denoted by $\MCV$.  The set of all
translation bounded measures is denoted by $\MTB$. As a subset of
$\MM$, it carries the vague topology. $\MCV$ is compact in this
topology.  There is an obvious action of $G$ on $\MTB$, again denoted
by $\alpha$, given by
\begin{equation*}
   \alpha \! : \; G\times  \MTB \; \longrightarrow \; \MTB 
   \, , \quad (t,\nu) \, \mapsto \, \alpha^{}_t\ts \nu 
   \quad \mbox{with} \quad (\alpha^{}_t \ts \nu)(\varphi) \, := \,
   \nu(\delta_{-t} \ast \varphi)
\end{equation*}
for $\varphi \in C_c (G)$.  Here, $\delta_t$ denotes the unit point
mass at $t \in G$ and the convolution $\omega \ast \varphi$ between
$\varphi\in C_c (G)$ and $\omega \in \MTB$ is defined by
\[ \omega \ast \varphi (s) := \int \varphi (s - u) \dd \omega(u).\]

It is not hard to see that $\alpha$ is continuous when restricted to a
compact subset of $\MTB$.

\begin{definition} 
  $\Oomega$ is called a dynamical system on the translation bounded 
  measures on\/ $G$  {\rm (TMDS)} if  
  $\varOmega$ is a compact  $\alpha$-invariant subset of\/ $\MCV$ for some open relatively compact $V$ and $C>0$.
 \end{definition}

Every translation bounded measure $\nu$ gives rise to a (TMDS)
$(\varOmega(\nu),\alpha)$, where
\[ \varOmega(\nu):=\overline{\{\alpha_t\nu : t\in G\}}.\]
More precisely, if $\nu\in \MCV$, then $\varOmega(\nu) \subset \MCV$.

\smallskip

As usual $\varphi \in C_b (G)$ is called {\em almost periodic} (in the
sense of Bohr) if, for every $\epsilon >0$, the set of $t\in G$ with
$\|\delta_t \ast \varphi - \varphi \|_\infty \leq \epsilon$ is
relatively dense in $G$. By standard reasoning this is equivalent to
$\{\delta_t \ast \varphi : t\in G\}$ being relatively compact in in
$C_b (G)$ (see e.g. \cite{Zai}). 

\begin{definition} 
A translation bounded measure $\nu$ is called strongly almost periodic
  if $\nu \ast \varphi$ is almost periodic (in the Bohr sense) for
  every $\varphi \in C_c (G)$.
\end{definition}

\subsection{Diffraction theory}\label{difftheo}

Having introduced our models, we can now discuss some key issues
of diffraction theory, where we follow \cite{BL,BL2,BM}. 

\smallskip

Let $\Oomega$ be a TMDS, equipped with an $\alpha$-invariant measure $m$. Fix
$\omega\in\Omega$ and let $\psi \in C_c (G)$ with $\int \psi(t) \dd t =1$ be given. Then, 
$\gamma_m : C_c (G) \longrightarrow \CC$ defined by
\[ \gamma_m (\varphi) := \int_\varOmega \int_G \int_G \varphi ( s + t) \psi
(t) \dd \omega(s) \dd \widetilde{\omega}(t) \dd m (\omega)\]
is a positive definite measure which does not depend on $\psi$ (provided
$\int \psi(t) \dd t =1$). Here, for $\nu\in \MM$, the measure
$\widetilde{\nu}$ is defined by $\widetilde{\nu}(\varphi):=\overline{ \nu(
  \overline{ \varphi(-\cdot)} }$. The measure $\gamma_m$ is called
\textit{autocorrelation measure}. Its Fourier transform exists and is called
\textit{diffraction measure} (see \cite{Fol, GdeL} for definition and background on Fourier transforms on measures). This measure describes the outcome of
actual diffraction experiments \cite{Cowley,Hof}. 
If  $\Oomega$ is ergodic,  $\gamma_m$ can be calculated via a limiting
procedure \cite{BL,Gouere-1,Hof,Martin2}.  
More precisely, recall that the convolution $\mu \ast \nu$ of two bounded 
measures $\mu$ and $\nu$ on $G$ is defined to be the measure $\mu\ast \nu 
(\varphi) :=\int \int \varphi(s+t) \dd\mu(s) \dd\nu(t)$. 
Now, if $\Oomega$ is uniquely ergodic 
\begin{equation}\label{gamma_als_mittel}
\gamma_m = \lim_{n\to \infty} \frac{1}{m_G(B_n)} \omega_{B_n} \ast
\widetilde{\omega_{B_n}}
\end{equation}
for every $\omega \in \Omega$ \cite[Thm.~5]{BL}. 
Here, the limit is taken in the vague topology,  $\omega_{B_n}$ denotes the
 restriction of $\omega$ to $B_n$ and $(B_n)$ is a van Hove sequence in $G$. 
This means \cite{Martin2} that for every  compact $K\subset G$, 
$$\lim_{n\to \infty} \frac{m_G (\partial^K B_n)}{m_G (B_n)} =0,$$
where for arbitrary $A,K\subset G$ we set
\begin{equation}\label{rand} \partial^K A:= (( K + A )\setminus A^\circ ) 
  \cup (( - K +\overline{G\setminus A})\cap A ),
\end{equation}
where the bar denotes the closure of a set and the circle denotes the
interior. If $\Oomega$ is uniquely ergodic,
we write $\gamma$ instead of $\gamma_m$. We also recall the following result from \cite{BL}. 

\medskip

\begin{theorem}\label{aequivalenz}  
  Let\/ $\Oomega$ be a\/ {\rm TMDS} with invariant measure\/ $m$. 
  Then, the following assertions are equivalent.
\begin{enumerate}
\item[(i)]  The measure\/ $\widehat{\gamma_m}$ is a pure point measure.
\item[(ii)] $T^\varOmega$ has pure point dynamical spectrum.  
\end{enumerate}
In this case, the group generated by $\{\lambda\in \Ghat: \gammahat(\{\lambda\})>0\}$ is the group of eigenvalues of $T^\varOmega$. 
\end{theorem}

\medskip

\section{Cut and project schemes for measures: Main results} \label{Cut}

In this section, we introduce cut and project schemes and discuss 
our main results. 

As usual a triple $\cp$ is called a cut-and-project scheme if 
$G$ and $H$  are locally compact $\sigma$-compact abelian groups and  
$\widetilde{L}$ is a {\em lattice\/} in $G\times H$ (i.e., a cocompact
discrete subgroup) such that 
\begin{itemize}
\item the canonical projection $\pi : G\times H \longrightarrow G$ is
one-to-one between $\widetilde{L}$ and $L:=\pi (G)$ (in other words,
$\widetilde{L}\cap (\{0\}\times H) =\{(0,0)\}$), and
\item the image $L^\star = \pi^{}_{\rm int}(\widetilde{L})$ of the canonical
  projection  $\pi^{}_{\rm int} : G\times H
\longrightarrow H$ is dense in $H$.
\end{itemize}
The group $H$ is called the {\em internal} space. Given 
these properties of the projections $\pi$ and $\pi^{}_{\rm int}$, one
can define the $\star$-map as $(.)^\star\!: L \longrightarrow H$
via $x^\star := \big( \pi^{}_{\rm int} \circ (\pi|_L)^{-1}\big) (x)$,
where $(\pi|_L)^{-1} (x) = \pi^{-1}(x)\cap\widetilde{L}$, for all $x\in L$. 
This situation can  be summarised in the following diagram.
\begin{equation*}
\begin{array}{cccccl}
    G & \xleftarrow{\,\;\;\pi\;\;\,} & G\times H & 
        \xrightarrow{\;\pi^{}_{\rm int}\;} & H & \\
   \cup & & \cup & & \cup & \hspace*{-2ex} \mbox{\small dense} \\
    L & \xleftarrow{\; 1-1 \;} & \widetilde{L} & 
        \xrightarrow{\,\;\quad\;\,} & L^\star & \\
   {\scriptstyle \parallel} & & & & {\scriptstyle \parallel} \\
    L & & \hspace*{-38pt}
    \xrightarrow{\hspace*{47pt}\star\hspace*{47pt}} 
    \hspace*{-38pt}& & L^\star
\end{array}
\end{equation*}
A cut and  project scheme gives rise to a dynamical system in the following
way: Define $\TT :=(G \times H) / \widetilde{L}$.  By assumption on $
\widetilde{L}$, $\TT$ is a compact abelian group. 
Let 
\[ G \times H \longrightarrow \TT, \;\:(t,k)\mapsto [t,k],\]
be the canonical quotient map.  There is a canonical continuous group
homomorphism
\[ \iota : G\longrightarrow \TT, \;\: t \mapsto [t,0].\] 
The homomorphism $\iota$ has dense range as $L^\star$ is dense in $H$. It
 induces an action $\beta$ of $G$ on $\TT$ via
\[\beta : G\times \TT \longrightarrow\TT,\:\; \beta_t([s,k]):= \iota(-t) +  [s,k]= [s - t,k].\]
The dynamical system $(\TT,\beta)$ will play a crucial role in our
considerations.  It is minimal and uniquely ergodic, as $\iota$ has
dense range.  Moreover, it has pure point spectrum. More precisely,
the dual group $\widehat{\TT}$ gives a set of eigenfunctions, which form
a complete orthonormal basis by Peter-Weyl theorem  (see \cite{Martin2} for
further details). Later we will also meet the canonical injective 
group homomorphism 
$$\kappa : H \longrightarrow \TT, \;\:h \mapsto [0,h].$$

\begin{definition}\label{def:admiss}
{\rm  (a)} A quadruple $\mcp$ is called a {\em measure cut and project scheme} 
if $\cp$ is a cut and project scheme and $\rho$ is an $\widetilde{L}$-invariant 
Borel measure  on $G\times H$. \\
{\rm (b)} Let $\mcp$ be a measure cut and project scheme. A function $f : H
\longrightarrow \CC$ is called {\em admissible} if it is measurable, locally 
bounded and for arbitrary $\varepsilon>0$ and $\varphi \in C_c (G)$ there 
exists a compact $Q\subset H$ with
\[ \int_{G\times H} |\varphi (t+s) f(h + k )| (1- 1_{Q} (h + k)) \dd |\rho|(t,h)\leq 
\varepsilon \]
for every $(s,k)\in G\times H$, where $1_{Q}$ denotes the characteristic function 
of $Q$. 
\end{definition}

An example of a measure cut and project scheme is given by a cut and project
scheme  $\cp$ and $\rho:= \delta_{\widetilde L} := \sum_{x\in\widetilde{L}} 
\delta_x$.  It is not 
hard to see that then every Riemann integrable $f : H\longrightarrow \CC$ is 
admissible. In this way all the ``usual'' cut and project schemes fall within 
measure cut and project schemes, see Section \ref{Dense} and Section 
\ref{Complementary} for details. 

\medskip

Our focus here will be to investigate admissible functions which are
continuous. However, some of our results will  hold for arbitrary admissible
functions.

\smallskip

Let  a measure cut and project scheme  $\mcp$  with an admissible 
$f$ be given. As shown in Proposition \ref{tb} below, the map 
\begin{equation}\label{def_nu}
\nu_f : C_c (G) \longrightarrow \CC, \:\;\varphi \mapsto \int_{G\times H} 
\varphi (t) f(h)
\dd \rho (t,h), 
\end{equation}
is a translation bounded measure.  Thus, we can consider its hull
\[\varOmega (\nu_f):=\overline{\{\alpha_t (\nu_f) : t\in G\}}.\]
By the discussion of the previous section, $(\varOmega (\nu_f),\alpha)$ is then
a TMDS. 

\medskip

Our main results  are the following three. The first deals with the
dynamical system side of the problem, the second and third deal with
diffraction. 

\begin{theorem} \label{main} Let  a measure cut and project scheme  $\mcp$
  with a continuous admissible $f$ be given and $\nu_f$ be defined as in
  \eqref{def_nu}. Then, the following assertions hold.

\begin{itemize}

\item[(a)] $\nu_f$ is strongly almost periodic.  In particular, $\varOmega
  (\nu_f)$ has a unique abelian group structure such that $G\longrightarrow\varOmega
  (\nu_f)$, $t\mapsto \alpha_t \nu_f$, is a continuous group homomorphism. 

\item[(b)] $(\varOmega (\nu_f),\alpha)$ is a factor of $(\TT,\beta)$ with
  factor map 
$ \mu : \TT \longrightarrow \varOmega(\nu_f)$ given by
  $\mu ([s,k]) (\varphi) =\int f(h+k) \varphi (s + t) \dd \rho(t,h)$. 
In fact, $\mu $ is a group homomorphism. 

\item[(c)] $(\varOmega (\nu_f),\alpha)$ is minimal, uniquely ergodic and has
  pure point spectrum with continuous eigenfunctions. The set of eigenvalues
  is contained in $\{\lambda\circ \iota : \lambda \in \widehat{\TT} \}\subset
  \widehat{G}$. 
\end{itemize}
\end{theorem}

Before we can state the next results, we recall the following result
on disintegration \cite[Sec.~33]{Loomis}. Let $\mcp$ be a measure cut 
and project scheme. Let for $\xi= [s,h] \in \TT$, the (well defined!) 
measure $\sigma_\xi$ on $G\times H$ be given by $\sigma_\xi (g) = 
\sum_{(l,l^\star)\in\widetilde{L}} g(s + l,h + l^\star)$ for $g\in 
C_c(G\times H)$. Then, there exists a unique measure $\rho_\TT$ on $\TT$
with
\begin{equation}\label{form:rt}
\int_{G\times H} g(s,h) \dd\rho(s,h) = \int_\TT \sigma_\xi (g)\dd\rho_\TT
(\xi)
\end{equation}
for all $g\in C_c(G\times H)$. 
In fact, (and this shows both existence and uniqueness) the measure
$\rho_\TT$ satisfies 
$$ \rho_\TT (b) =\int_{G\times H}  b([s,h]) \chi_{Z} (s,h) \dd\rho(s,h)$$
for $b\in C(\TT)$, whenever $Z$ is a fundamental cell of $\widetilde{L}$ 
in $G\times H$ (i.e., $Z$ is a measurable subset of $G\times H$ such that
$Z\longrightarrow \TT$, $(s,h)\mapsto[s,h]$, is bijective.)

Moreover, for a function $f:H\to\mathbb C$ define
$\widetilde{f}(h):=\overline{f(-h)}$. For a measure $\rho$ on $G\times
H$, define $\overline{\rho}$ by $\overline{\rho} (g) = \overline{ \rho
(\overline{g})}$ for every $g\in C_c(G\times H)$.

\begin{theorem}\label{main_autocorrelation}  
Let a measure cut and project scheme $\mcp$ with a continuous
admissible $f$ be given. Then, $f$ is integrable and in particular
$(f\ast \widetilde{f}) (h)$  exists for almost every $h$ in $H$.  For
$\varphi\in C_c (G)$ and $\xi= [s,h]\in\TT$ define
$$\gamma_\xi (\varphi) :=\frac{1}{(m_G\times m_H)(Z)}\sum_{(l,l^\star) \in \widetilde{L}}
(f\ast\widetilde{f})(h - l^\star) \varphi (s - l)$$ 
whenever this
exists. Note that this is well defined, i.e., $\gamma_\xi (\varphi)$ 
does not depend on the chosen representative of $\xi$.  Then, for every $\varphi \in
C_c (G)$ and $\rho \times \overline{\rho}$ almost every $(\xi,\eta)$,
$\gamma_{\xi-\eta} (\varphi)$ exists and the autocorrelation $\gamma$
of $\nu_f$ satisfies
$$ \gamma(\varphi) = \int \int \gamma_{\xi-\eta} (\varphi) \dd \rho_\TT
(\xi) \dd\overline{\rho_\TT} (\eta).$$ 
If $\rho =\delta_{\widetilde{L}}$, then $f$ is square integrable, $(f\ast
\widetilde{f}) (h)$ exists for every $h\in H$ and $\gamma = \gamma_0$.
\end{theorem}

\begin{theorem}\label{main_diffraction} 
Let a measure cut and project scheme $\mcp$ with a continuous 
admissible $f$ be given and $\nu_f$ be defined as in \eqref{def_nu}.  
Let $\gammahat$ be the associated diffraction measure. Then, 
$\gammahat$ is a pure point measure supported on 
$\{\lambda\circ \iota : \lambda \in \widehat{\TT} \}$. More precisely,
$\gammahat= \sum_{\lambda\in \widehat{\TT}} |c_\lambda|^2
\delta_{\lambda\circ \iota}$
with
$$ c_\lambda := \frac{\rho_\TT ( \lambda) }{(m_G \times m_H)_\TT (1)}
\int_H f(h) (\lambda\circ \kappa) (h) {\rm d}h = \lim_{n\to\infty} 
\frac{1}{m_G (B_n)} \nu_f (\chi_{B_n} \cdot\overline{(\lambda\circ \iota)})$$
for every van Hove sequence $(B_n)$.
\end{theorem}


\section{Strongly almost periodic measures} \label{Strongly} 

In this section, we show that almost periodic measures on $G$ give
rise to topological groups, which are dynamical systems. Much of the
material of this section can is at least implicitly contained in the
literature, particularly in \cite{GdeL}. However, for the
convenience of the reader, and as our dynamical system perspective is
not the usual approach to these results, we include proofs.

\smallskip

We start by introducing the relevant topology.  For $\nu\in \MTB$ and
$\varphi \in C_c (G)$, the convolution $\nu \ast \varphi$ belongs to
$C_b (G)$.  Let $C_b (G)$ be equipped with the supremum norm.  We then
define the strong topology $\cT_s$ to be the weakest topology on
$\MTB$ such that all maps
$$
\MTB\longrightarrow C_b (G), \;\:\nu\mapsto \nu\ast \varphi,  $$
are
continuous.  Alternatively, we can describe this topology by considering
$\MTB$ as a subset of $\prod_{\varphi \in C_c (G)} C_b (G)$ which is equipped
with the product topology via $$
i : \MTB \longrightarrow \prod_{\varphi \in
C_c (G)} C_b (G), \;\: i(\nu) := (\varphi \mapsto \nu \ast \varphi).$$
The projection on the $\psi$ component 
$$\prod_{\varphi \in C_c (G)} C_b (G)\longrightarrow C_b (G), \;\:  x \mapsto x_\psi,$$
is denoted by $p_\psi$.

\begin{prop} $i(\MTB)$ is closed in $\prod_{\varphi \in C_c (G)} C_b (G)$.
\end{prop} 

\begin{proof} Let  $(\nu_n)$ be a net in $\MTB$ such that $i (\nu_n)$
  converges to $x\in \prod_{\varphi \in C_c (G)} C_b (G)$.  Then, $\nu_n$
  converge in particular in the vague topology to a measure $\nu$ and $\nu
  \ast \varphi$ belongs to $C_b (G)$ for every $\varphi \in C_c (G)$.  Thus,
  $\nu$ is translation bounded and it is not hard to see that $i(\nu) =x$.
\end{proof}

\begin{lemma}\label{lemma_sap} 
Let $\nu \in \MTB$ be given.  The following assertions are equivalent:
\begin{itemize}
\item[(i)]  The measure $\nu$ is strongly almost periodic (i.e. $\nu \ast
  \varphi$ is Bohr almost periodic for every $\varphi \in C_c (G)$). 

\item[(ii)] $\{\alpha_t\nu : t\in G\}$ is relatively compact in $\cT_s$. 

\item[(iii)] The topological space $\varOmega(\nu)$ (the hull of $\nu$ in the
  vague topology)  is a topological group
  with addition $\dotplus$ satisfying $\alpha_s \nu \dotplus \alpha_t \nu =
  \alpha_{s+t} \nu$ for all $ s,t\in G$. 
\end{itemize}

\end{lemma}
\begin{proof} 

(ii) $\Longrightarrow $ (i): Define $B := \{ \alpha_t \nu : t\in G\}$
and $A_\varphi :=\{ \delta_t \ast (\nu \ast \varphi): t \in G\}$ and
$C_\varphi :=\overline{A_\varphi}$. Then, $i(B)$ is relatively compact
by (ii) and the definition of the strong topology.  A direct calculation
shows $A_\varphi = p_\varphi (i (B))$.  As $p_\varphi$ is continuous,
compactness of $C_\varphi$ follows now from $C_\varphi= \overline{
A_\varphi} = \overline{ p_\varphi (i(B))} = p_\varphi (
\overline{i(B)})$.

\smallskip


(i) $\Longrightarrow$ (iii):  For $R>0$, let $K(R)\subset \CC$
  be the closed ball around the origin with radius $R$. As $\varOmega(\nu)$ is
  compact, there exists for each $\varphi \in C_c (G)$ an $R_\varphi>0$ with
\[ \omega (\varphi) \in K(R_\varphi)\]
for every $\omega \in \varOmega(\nu)$. We can and will therefore consider
$\varOmega(\nu)$ to be a compact subset of $\prod_{\varphi \in C_c (G)}
K(R_\varphi)$.  For $\varphi \in C_c (G)$ define $\check{\varphi}\in C_c (G)$
by $\check{\varphi} (t) := \varphi (-t)$.  As $\nu$ is strongly almost periodic, for each $\varphi \in C_c
(G)$, the function $v_\varphi := \nu \ast \check{\varphi}$
is almost periodic.  Thus, the closure $\varOmega_\varphi$ of $\{\delta_t \ast
v_\varphi : t\in G\}$  with respect to the supremum norm $\|\cdot\|_\infty$
is a compact  abelian group (see e.g. \cite{HR}) . Moreover,  $j_\varphi : G\longrightarrow
\varOmega_\varphi$, $t\mapsto \delta_t \ast v_\varphi $, is a continuous group
homomorphism.  Then, $\prod_{\varphi \in C_c (G) } \varOmega_\varphi$ is an abelian group, which is compact by Tychonov's theorem. Moreover,   
\[ j : G \longrightarrow \prod_{\varphi \in C_c (G) } \varOmega_\varphi, \;\:
t\mapsto (\varphi \mapsto j_\varphi (t)),\]
is a continuous group homomorphism.  Obviously, $j(G)$ is a subgroup of
$\prod_{\varphi \in C_c (G) } \varOmega_\varphi$.  In particular, its closure
$\overline{j(G)}$ is a compact abelian group. We show that $\varOmega(\nu)$ is
homeomorphic to $\overline{j(G)}$:

Consider the evaluation at $0$ 
\[ \Delta : \prod_{\varphi \in C_c (G) } \varOmega_\varphi\longrightarrow
\prod_{\varphi \in C_c (G)} K(R_\varphi), \;\:(\varphi \mapsto w_\varphi)
\mapsto (\varphi \mapsto w_\varphi (0)).\]
Then, $\Delta$ is continuous, as each $\varOmega_\varphi$ is equipped with the
supremum norm. A short calculation shows
\[ \Delta \circ j (t) = (\varphi \mapsto v_\varphi (-t)) = (\varphi \mapsto
(\alpha_t \nu) (\varphi)) = \alpha_t \nu.\]
As $\Delta$ is continuous and $\overline{j(G)}$ compact, this gives
\[ \Delta (\overline{ j(G)}) = \overline{\Delta (j(G))} = \overline{ \{ \alpha_t
\nu : t\in G \}   } = \varOmega (\nu). \]
Thus, $\Delta$  maps $\overline{j(G)}$ onto
$\varOmega(\nu)$. 

\smallskip

\textbf{Claim.} $\Delta$ is one-to-one on $\overline{j(G)}$: 

Proof of claim.  As every $w=(\varphi \mapsto w_\varphi)\in j(G)$ satisfies  $w_\varphi (t) = w_{\varphi(\cdot -t)} (0)$, the same holds for $w\in \overline{j(G)}$. Thus, the evaluations at $0$ determine all coordinates and $\Delta$ is injective on $\overline{j(G)}$. 

\smallskip

These considerations show that $\Delta : \overline{j(G)}
\longrightarrow \varOmega(\nu)$ is a homeomorphism.  As
$\overline{j(G)}$ is a compact abelian group, $\varOmega(\nu)$
inherits the structure of a compact abelian group as well, and
\[ \alpha_t \nu \dotplus \alpha_s \nu = \Delta (j(s)) \dotplus \Delta(j(t)) =
\Delta(j(s) \dotplus j(t)) = \Delta (j(s+t)) = \alpha_{t+s} \nu.\]
This finishes the proof of this implication.

\smallskip

(iii) $\Longrightarrow$ (ii) : As in the proof of (ii)$\Longrightarrow$ (i), we
define $B := \{ \alpha_t \nu : t\in G\}$ and note that relative compactness of
$B$ is equivalent to relative compactness of $i(B)$. 

Choose an arbitrary $\varphi \in C_c (G)$. 
We have to show that the closure $\varOmega_\varphi$ of 
\[ \{\delta_s \ast \nu \ast \varphi : s\in G\}  \]
in $(C_b (G), \|\cdot\|_\infty)$ 
is compact. Then,  $\prod_{\varphi \in C_c (G)} \varOmega_\varphi$ is compact
by Tychonov's theorem, and relative compactness of $i(B) \subset \prod_{\varphi
\in C_c (G)} \varOmega_\varphi$ follows. 

Obviously, the function $\alpha_\varphi : \varOmega \longrightarrow
\CC$, $\omega \mapsto \omega \ast \varphi (0)$, is continuous. As $\varOmega$
is a compact group, the map
\[ \varOmega \longrightarrow C(\varOmega), \omega\mapsto \alpha_\varphi(\omega
\dotplus\cdot)\]
is then continuous. Moreover, 
\[C(\varOmega)\longrightarrow C_b (G), b \mapsto
(t \mapsto b(\alpha_t\nu)),\]
is continuous. Using $\alpha_t\nu \dotplus  \omega = \alpha_t \omega$, we
infer that 
\[ p_\varphi : \varOmega \longrightarrow C_b (G), \, \omega \mapsto (t\mapsto
\alpha_\varphi (\alpha_t \omega)),\]
is continuous as composition of continuous maps. In particular, $p_\varphi
(\varOmega)$ is compact. A direct calculation shows $p_\varphi (\alpha_s \nu)
= \delta_s \ast \nu \ast \varphi$ yielding  $\varOmega_\varphi \subset
p_\varphi (\varOmega)$, and compactness of $\varOmega_\varphi$ follows. 
\end{proof}

\section{Factors}\label{Factors}

In this paper we study properties of $\varOmega (\nu)$ for suitable 
translation bounded measures $\nu$. Existence of a factor map 
\begin{equation}\label{existence-mu} \mu : \TT \longrightarrow \varOmega (\nu),
\end{equation}
with a suitable abelian compact group $\TT$ will be of key importance. In this
section, we provide abstract background to existence and use of such a
factor map.

\smallskip

Factors inherit basic features of their underlying dynamical systems. The
following statements summarises results proved in Section 3 of \cite{BL2}.

\begin{prop}\cite[Sec.~3]{BL2} \label{thm_factor} 
  Let\/ $\Oomega$ be a topological dynamical system and let\/ $\Ttheta$ be a
  factor with factor map $\varPhi:\Omega\to\Theta$. For an invariant probability 
  measure $m$ on $\Oomega$, we define the invariant probability measure
  $\varPhi(m)$ on $\Theta$ by $\varPhi (m) (g) := m( g\circ \varPhi)$ for $g\in
  C_c(\varTheta)$. The following assertions hold.
\begin{itemize}
\item[\rm (a)] If $\Oomega$ is uniquely ergodic with invariant probability measure $m$, then $\Ttheta$ is uniquely ergodic with unique invariant
probability measure $\varPhi(m)$.
\item[\rm (b)] If\/ $\Oomega$ has pure point dynamical spectrum when equipped
  with the invariant probability measure $m$, $\Ttheta$ has pure point
  dynamical spectrum when equipped with the measure $\varPhi (m)$.
\item[\rm (c)] If\/ $\Oomega$ is minimal, so is\/ $\Ttheta$.
\item[\rm (d)] If\/ $\Oomega$ is uniquely ergodic with pure point 
   dynamical spectrum and all of its eigenfunctions continuous, 
   the same holds for\/ $\Ttheta$. 
\end{itemize}
\end{prop}

We will be interested in special factors of compact groups. The relevant lemma is the following. 

\begin{lemma}\label{groupstructure}  
  Let a compact abelian group $\TT$ and a continuous group
  homomorphism $\iota :G\longrightarrow \TT$ with dense range be
  given. Let $\beta$ be the associated action of $G$ on $\TT$
  i.e. $\beta_t (\xi)=\iota(t)\xi$ for $\xi \in \TT$ and $t\in G$.  If
  $\Oomega$ is a factor of $(\TT,\beta)$ with factor map $\mu$, then
  there exists a unique topological group structure on $\varOmega$
  such that $\alpha_t\mu(0) \dotplus \alpha_s \mu(0) = \alpha_{t
  +s}\mu (0)$. With respect to this group structure $\mu$ is a group homomorphism.
\end{lemma}
\begin{proof} By denseness, there can be at most one group structure with 
\begin{equation*}\label{orbit}
\alpha_t\mu(0)
  \dotplus \alpha_s \mu(0) = \alpha_{t +s}\mu (0).
\end{equation*}
Next, we show that
\begin{equation}\label{subgroup}
\mu(\eta) = \mu (\rho) \Longleftrightarrow \mu (\rho - \eta) =\mu(0). 
\end{equation}
 Let $\eta =\lim \beta_{t_s} (0)$. 

``$\Rightarrow$:'' We have $- \eta =
\lim \beta_{-t_s} (0)$ and
\[\mu(\rho - \eta) = \lim \mu (\beta_{- t_s} \rho) = \lim \alpha_{-t_s} \mu (\rho) = \lim \alpha_{-t_s} \mu (\eta) = \lim \mu (\beta_{-t_s} \eta) =\mu (0). \]

\smallskip

``$\Leftarrow$:'' Set $\xi :=\rho -\eta$. Thus, $\xi + \eta =\rho$ and 
 \[\mu (\rho) = \lim \mu ( \xi +  \iota(t_s)(0)) =\lim \mu ( \beta_{t_s} \xi)\\ = \lim \alpha_{t_s} \mu(\xi) =  \lim \alpha_{t_s} \mu(0)\\ = \lim \mu( \beta_{t_i} 0) = \mu (\eta).
\]

By continuity of $\mu$, the set
\[ U :=\{\xi \in\TT : \mu(\xi) = \mu(0)\}\]
is closed. Moreover, by \eqref{subgroup}, $U$ is a subgroup. As $G$ is abelian, $U$ is normal and 
\[\widetilde{\mu} : \TT/U \longrightarrow \varOmega, \, \rho + U \mapsto \mu(\rho), \]
is well defined, continuous and bijective. Thus, $\varOmega$ is homeomorphic to the group $\TT/U$ and inherits the desired group structure. 
\end{proof}

\begin{theorem}\label{SAP} 
Let a compact abelian group $\TT$ and a continuous group homomorphism
  $\iota : G\longrightarrow \TT$ with dense range be given. Let
  $\beta$ be the associated action of $G$ on $\TT$, i.e., $\beta_t
  (\xi)=\iota(t)\xi$. Let $\nu$ be a translation bounded measure on
  $G$. Then, the following assertions are equivalent:
\begin{itemize}
\item[(i)] There exists a factor map $\mu : \TT \longrightarrow \varOmega
  (\nu)$. 
\item[(ii)] The measure $\nu$ is strongly almost periodic,
  $\varOmega(\nu)$ is a topological group satisfying $\alpha_s\nu
  \dotplus \alpha_t \nu = \alpha_{s+t} \nu$ and
  $\widehat{\varOmega(\nu)} \subset \widehat{\TT}$, where both groups
  are considered as subgroups of $\Ghat$.
\end{itemize}
If there exists a cut and project scheme $(G,H,\widetilde{L})$ with $\TT=
(G\times H)/ \widetilde{L}$, this is equivalent to

\begin{itemize}

\item[(iii)] There exists an $\widetilde{L}$-invariant measure $\sigma$ on
  $G\times H$ together with a continuous disintegration $\sigma =\int_H
  \sigma_h  \, \dd h$ (i.e. $H\longrightarrow \MM$, $h\mapsto \sigma_h$, is
  continuous and $\int_{G\times H} g(s,h) \dd\sigma(s,h) = \int_H \sigma_h
  (g(\cdot,h)) \dd h$ for every $g\in C_c (G\times H)$) and $\nu = \sigma_0$. 

\end{itemize}

\end{theorem}

\begin{proof} 
(i)$\Longrightarrow$ (ii): As $\mu$ is a factor map, $\varOmega(\nu)$ 
inherits by Lemma \ref{groupstructure} the structure of an abelian group 
such that $\mu$ is a group homomorphism.  As $\mu$ is onto, the inclusion
  $\widehat{\varOmega(\nu)} \subset \widehat{\TT}$ follows. As  
$\varOmega(\nu)$  is a group, the measure $\nu$ is almost periodic by 
Lemma \ref{lemma_sap}.

\smallskip

(ii)$\Longrightarrow$(i): Dualising the inclusion
$\widehat{\varOmega(\nu)}\subset \widehat{\TT}$ gives a map
$\widehat{\widehat{\TT} } \longrightarrow \widehat{
\widehat{\varOmega(\nu)}}$. This map is onto as
$\widehat{\varOmega(\nu)}$ is closed as a subgroup of the discrete
group $\widehat{\TT}$ (see e.g. Corollary 4.41 in \cite{Fol}). By
Pontryagin duality, we obtain a surjective map from $\TT$ to
$\varOmega(\nu)$. It is easy to see that it is a factor map.

 (i) $\Longrightarrow$ (iii): For $h\in H$ define $\sigma_h
  :=\mu([0,h])$. Then, $\sigma_0=\nu$ by definition and $h\mapsto
  \sigma_h$ is continuous as $\mu$ is continuous. Moreover, $\sigma
  =\int_H \sigma_h \dd h$ is $\widetilde{L}$-invariant as $\alpha_l
  (\sigma_{h + l^\star}) = \mu([l,h+l^\star]) = \mu([0,h]) =\sigma_h$
  for every $(l,l^\star)\in \widetilde{L}$.

\smallskip

 (iii) $\Longrightarrow$ (i): Define $\mu' : G\times H \longrightarrow
 \MM$ by $\mu' (t,h) :=\alpha_t \sigma_h$. Then, $\mu'$ is continuous
 as $h\mapsto \sigma_h$ is and $\mu' (0,0)=\nu$. Moreover,
 $\widetilde{L}$-invariance of $\sigma$ gives
$$ \int_H \sigma_h (g(\cdot,h)) \dd h = \sigma (g) = \sigma(g (\cdot -l,
\cdot - l^\star)) = \int_H \alpha_l(\sigma_{h+l^\star})(g(\cdot,h)) \dd h$$
for every $g\in C_c (G\times H)$ and every $(l,l^\star)\in
\widetilde{L}$. As $h \mapsto \sigma_h$ is continuous, this shows
$\sigma_h = \alpha_l(\sigma_{h+l^\star})$ for every $(l,l^\star)\in
\widetilde{L}$, and $\widetilde{L}$-invariance of $\mu'$ follows.
 \end{proof}

\noindent\textbf{Example.} 
We show that the previous result covers the ``classic'' 
quasiperiodic functions: Let a cut and project
scheme $\cp$ be given with the associated canonical projection $p : G\times H
\longrightarrow \TT$. Let $ A : \TT \longrightarrow \RR$ be continuous and set
\[\rho := A\circ p\, \dd t\times\dd h.\]
Thus, $a : G \longrightarrow \RR$, $a(t):= A([t,0])$ is quasiperiodic. 

If we now define $\sigma_h := A([\cdot,h]) \dd t$, then $\sigma_0 = a \dd t$ and
$h\mapsto \sigma_h$ is continuous with $\rho = \int_H \sigma_h \dd h$. 
Thus, the previous theorem applies to $\nu :=
a \, \dd t$, and we obtain a factor map and group homomorphism $\mu$ satisfying
\eqref{existence-mu}.  In this sense, our results cover not only Dirac combs
but also quasiperiodic functions. Of course, for $a$ quasiperiodic validity
of \eqref{existence-mu} can also rather directly be shown from the definitions.

\section{A study of admissibility} \label{Study}

The aim of this section is to study admissibility. 

\begin{prop} \label{rho_tb} 
Let $\mcp$ be a measure cut and project scheme. Then, the  measure 
$\rho$ is translation bounded.
\end{prop}

\begin{proof}  Let $C$ be the closure of a nonempty open relatively 
compact set in $G\times H$. Let $D$ be a compact subset of $G\times H$ 
containing a fundamental domain of $\widetilde{L}$. Then, by  
$\widetilde{L}$-invariance of $\rho$, for arbitrary $u\in G\times H$, 
there exists an $v\in D$  with $ |\rho| (u + C) = |\rho| (v + C)$. As 
$D + C$ is compact and $\rho$ is a Borel measure, this implies
\[ |\rho| (u + C) \leq  |\rho| (D + C) =  const<\infty\]
for every $u\in G\times H$. The proposition follows.
\end{proof}

\begin{prop} 
Let $\mcp$ be a measure cut and project scheme and $f : H\longrightarrow 
\mathbb C$ be locally bounded and measurable. Then, the following 
assertions are equivalent:

\begin{itemize}
\item[(i)] The function $f$ is admissible. 
\item[(ii)]  For all $\varepsilon>0$ and $\varphi  \in C_c
  (G)$ there exists a $\chi : H \longrightarrow [0,1]$ in $C_c (H)$  with
\[ \int_{G\times H} |\varphi (t+s) f(h+k)|( 1- \chi (h+k)) | 
\dd |\rho|(t,h)\leq \varepsilon \]
for every $s\in G$ and $k\in H$.
\item[(iii)] For all $\varepsilon>0$ and $\varphi \in C_c
  (G)$ there exists a $g\in C_c (H)$  with
\[ \int_{G\times H} |\varphi (t+s) (f(h+k) - g(h+k)) |  \dd |\rho|(t,h)
\leq \varepsilon \]
for every $s\in G$ and $k\in H$.
\end{itemize}
\end{prop}

\begin{proof} 
(i)$\Longrightarrow $ (ii): By Tietze's extension theorem, there exists a 
 $\chi\in C_c (H)$  with $1\geq \chi \geq 0$ and $\chi=1$ on $Q$. Then  $1 -
 1_Q \geq 1 - \chi \geq 0$.\\
(ii)$\Longrightarrow $ (iii): Set  $g= \chi f$.\\
 (iii) $\Longrightarrow $ (i): Set $Q:=\supp(g)$. 
\end{proof}

\medskip 

\begin{prop}\label{tb} 
Let $\mcp$ be a measure cut and project scheme and $f : H\longrightarrow 
\mathbb C$ be admissible. For every $k\in H$ and $\varphi 
\in C_c (G)$, there exists a $C_\varphi\geq 0$ with 
$$ 
\int |\varphi (s + t) f(h + k)| \,{\rm d}|\rho|(t,h)\leq C_\varphi 
$$
for every $ s \in G$. 
In particular, for every $(s,k)\in G\times H$ the map 
\[ \mu'(s,k) : C_c (G) \longrightarrow \CC, \;\: \varphi \mapsto \int_{G\times H} f(h + k)
\varphi (t + s) \dd \rho (t,h) \]
is a translation bounded measure on $G$ and so  is $\nu_f = \mu'(0,0)$. 
\end{prop}

\begin{proof} Choose $k\in H$ and let $\varphi$ in $C_c (G)$ be given. 
As  $f$ is admissible, we can find a continuous  $\chi :
H\longrightarrow [0,1]$ with compact support   with $ \int_{G\times H}
|\varphi (s + t)  f(k + h)|   (1-\chi (h +k))   \dd |\rho|(t,h)\leq 1$ for
every $s\in G$. 
Moreover, as $\rho$ is translation bounded   and $(t,h)\mapsto |\varphi (s+t)
f(h+k) \chi (h+k)|$ is  bounded  with compact support, there exists a $
C'$ with
$$ \int_{G\times H} |\varphi (s+t) f(h+k) \chi (h+k)| \dd
|\rho|(t,h)\leq C'.$$ The first statement follows with $C:=1 + C'$.
Translation boundedness of the measures $\mu'(s,k)$ now follows easily
by choosing nonnegative $\varphi$ which are equal to $1$ on the closure
of an arbitrary open relatively compact $V$.
\end{proof}

\begin{prop}\label{admissible-fat}
Let $\mcp$ be a measure cut and project scheme and $f : H\longrightarrow \mathbb C$ 
be admissible. For every $K\subset H$ compact, $\varepsilon >0$ and $\varphi \in
  C_c (G)$, there exists a compact $Q_K$ with
\[ \int_{G\times H} |f(h + k) \varphi (t+s) (1 -  1_{Q_K}(h))| \dd |\rho|(t,h)\leq
\varepsilon \] 
for every $s\in G$ and $k\in K$. 
\end{prop}

\begin{proof} 
As $f$ is admissible, we can find  $\chi : H \longrightarrow [0,1]$ 
continuous with compact support with
\[ \int_{G\times H}  | \varphi(t + s) f(h+k) (1- \chi(h+k)) | 
\dd |\rho|(t,h)\leq \varepsilon\]
for every $s \in G$ and $k\in H$.  Then, $Q_K := \supp(\chi) -K$ has 
the desired properties.
\end{proof}

So, far our discussion of admissibility did not assume continuity of $f$. We
will now come to a characterisation of admissibility for continuous
$f$. 

\begin{prop}\label{Fatou} 
Let $\mcp$ be a measure cut and project scheme with
  metrisable $H$.  Let $c> 0$ and $ g : G\times H \longrightarrow [0,\infty)$ 
continuous be given with
\[ \int_{G\times H} g (t + s,h)  \dd |\rho|(t,h)\leq c \]
for every $s\in G$. Then, 
\[ \int_{G\times H} g (t + s,h+k)  \dd |\rho|(t,h)\leq c \]
for every $(s,k)\in G\times H$.
\end{prop}

\begin{proof} As  $\rho$ is $\widetilde{L}$-invariant, the assumption implies 
\[ \int_{G\times H} g(t+s + l, h + l^\star) \dd |\rho|(t,h)\leq c \]
for every $s\in G$ and every $(l,l^\star) \in \widetilde{L}$.  As $s$ is arbitrary,
this implies 
\[ \int_{G\times H} g(t+s, h + l^\star) \dd |\rho|(t,h)\leq c \]
for every $s\in G$ and every $l^\star \in \widetilde{L}$.  As $L^\star$ is dense
in $H$ by definition of a cut and project scheme and $H$ is metrisable by
assumption, we can find, for any $k\in H$, a sequence $(l_n^\star)$ in
$L^\star$ with $l_n^\star\to k$. Now, the statement follows from Fatou's Lemma 
by continuity of $g$. 
\end{proof}

\begin{prop}\label{char_admissible}
Let $\mcp$ be a measure cut and project scheme with metrisable $H$. Let 
$f: H\longrightarrow \CC$ be continuous. Then, $f$  
is admissible if and only if for arbitrary $\varepsilon>0$ and $\varphi 
\in C_c(G)$ there exists a compact $Q\subset H$ with
\[ \int_{G\times H} |\varphi (t+s) f(h  )| \,(1- 1_{Q} (h)) \dd |\rho|(t,h)\leq \varepsilon \]
for every $s\in G$.
\end{prop}

\begin{proof} 
The ``only if'' part is immediate from the definition of admissibility. 
The ``if'' part follows from the previous proposition.
\end{proof}

We finish this section by discussing restrictions on $f$ imposed by
the admissibility requirement.

\begin{prop} \label{integrability} 
Let $\mcp$ be a measure cut and project scheme and $f :
  H\longrightarrow \mathbb C$  admissible. Then, $f$ is integrable
  with respect to the Haar measure on $H$.
\end{prop}
\begin{proof} As $\rho$ is $\widetilde{L}$-invariant, we can choose continuous
  $ \varphi : G \longrightarrow [0,\infty)$ and $ \psi : H \longrightarrow
  [0,\infty)$ with compact support such that
\begin{equation}\label{lowerbound} \int \varphi ( s + t) \psi (k - h)\,
  {\rm d}|\rho| (t,h)\geq 1
\end{equation}
for all $(s,k)\in G\times H$.  

\smallskip

Set $C:= \int \psi \dd h$. By admissibility, there exists a compact
$Q\subset H$ such that
$$ \int \varphi (t + s) |f(h +k)| (1 - \chi_Q (h+k)) \,{\rm
d}|\rho|(t,h)\leq 1$$ 
for all $(s,k)\in G\times H$. This gives
$$ \int_H \psi (k)  \left(\int \varphi (t + s) |f(h +k)| (1 - \chi_Q (h+k))
\, {\rm d}|\rho|(t,h)\right) {\rm d}k\leq C.$$
Fubini's theorem and the translation invariance of the Haar measure then imply
$$ \int_H |f(k)| ( 1 - \chi_Q (k)) \left(\int \varphi (s + t) \psi (k-h ) \,{\rm d}|\rho|(
t,h)\right) {\rm d}k \leq C.$$
Now, the lower bound \eqref{lowerbound} implies
$$ \int_H |f(k)| ( 1 - \chi_Q (k)) \,{\rm d}k \leq C.$$
Thus, $ f(k) ( 1 - \chi_Q (k))$ is integrable. As $f$ is locally bounded, $f
\chi_Q$ is also integrable, and integrability of $f$ follows. 
\end{proof}

\noindent \textbf{Remark.} The converse of this proposition does not
hold, as can be seen by choosing $\rho:=\sum_{x\in \widetilde{L}}
\delta_x$ and and $f$ continuous and integrable with ``peaks''. More
precisely, the following holds.

\begin{prop} \label{squareintegrability}
Let $\mcp$ be a measure cut and project scheme with $\rho =
\delta_{\widetilde{L}}$. Let $f$ be continuous and admissible. Then,
$f$ is bounded and  square integrable. 
\end{prop}
\begin{proof} 
Let a nonnegative $\varphi \in C_c (G)$ with $\varphi(0)=1$ be
given. As shown in Proposition \ref{tb}, there exists $C_\varphi>0$
with
$$ \sum_{(l,l^\star)\in \widetilde{L}} \varphi(s +l) |f(l^\star)| \leq
C_\varphi$$
for every $s\in G$. As $L^\star$ is dense in $H$ and $f$ is
continuous, we infer
$$ \sup\{|f(h)| : h \in H\} =\sup\{|f(l^\star)| : l\in L\} \leq C_\varphi.$$
This proves boundedness of $f$. As $f$ is integrable by the previous
proposition, square integrability follows from boundedness.   
\end{proof}

In some sense, admissibility with respect to $\delta_{\widetilde{L}}$ 
implies admissibility with respect to any other $\widetilde{L}$-invariant measure. 

\begin{prop} \label{extension-of-admissibility}
Let $\mcp$ be a measure cut and project scheme. Let $f$ be continuous 
and admissible with respect to $(G,H,\widetilde{L}, 
\delta_{\widetilde{L}})$. Then, $f$ is admissible with respect to $\mcp$. 
\end{prop}
\begin{proof} To simplify the notation, we write $\zeta$ for 
$\delta_{\widetilde{L}}$. As $\rho$ is $\widetilde{L}$-invariant, there 
exists a finite measure $\rho_0$ supported on a fundamental domain $Z$ 
of $\widetilde{L}$ with 
$$ \rho =  \zeta \ast \rho_0.$$
As $f$ is admissible with respect to  $\zeta$, for each $\varphi \in 
C_c (G)$ and $\varepsilon >0$, there exists $Q\subset H$ compact with
$$ 
\int_{G\times H} |\varphi (t+s) f(h + k )| (1- 1_{Q} (h + k)) \dd 
|\zeta|(t,h)\leq \varepsilon
$$
for every $s\in G$ and $k\in H$. This gives 
\begin{displaymath}
\begin{split}
\int_{G\times H} |\varphi (t+s) &f(h + k )| (1- 1_{Q} (h + k)) 
\dd |\rho|(t,h) \leq \\
&\int_{Z}  \int_{G\times H}  |\varphi (t+s) f(h + k )| (1- 1_{Q} 
(h + k)) \dd |\zeta|(t,h) \dd|\rho_0 | (s,k) \leq |\rho_0| (Z) \,  \varepsilon
\end{split}
\end{displaymath}
for every $s\in G$ and $k\in H$, and we obtain admissibility of 
$f$ with respect to $\rho$. 
\end{proof}

\section{ Proof of Theorem \ref{main}} \label{Proof_main}

In this section we provide a proof of Theorem \ref{main}. 

\begin{lemma}\label{property_mu} 
Let $\mcp$ be a measure cut and project scheme and $f : H\longrightarrow 
\mathbb C$ be admissible and continuous. The map $\mu' : G\times H \longrightarrow
  \MTB$
\[ \mu'(s,k) : C_c (G) \longrightarrow \CC, \;\: \varphi \mapsto \int_{G\times H} f(h + k)
\varphi (t + s) \dd \rho (t,h) \]
defined in Proposition \ref{tb}  is continuous and
  $\widetilde{L}$-invariant, i.e., one has
  $\mu'(s+k,h+k^\star)=\mu'(s,h)$ for arbitrary $(s,h)\in G\times H$ and arbitrary
  $(k,k^\star)\in \widetilde L$.
\end{lemma}

\begin{proof} 
  Invariance is immediate from the definitions. Continuity follows easily from
  an $\varepsilon / 3$-argument: Namely, let $\{(s_n,k_n)\}$ be a net for which $
  (s_n,k_n) \longrightarrow (s,k)$ in $G\times H$.  We have to show
\[ \mu'(s_n,k_n) (\varphi) \longrightarrow \mu'(s,k) (\varphi)\]
for every $\varphi \in C_c (G)$.  Let $\varphi \in C_c (G)$ and $\epsilon>0$
be given.  As $\{k_n\}$ converges to $k$, there exists a compact neighbourhood
$K$ of $k$ with $\{k_n : n\geq n_0\}\subset K$. Therefore, by Proposition
\ref{admissible-fat}, there exists a compact $Q_K$ in $H$ with
\[ \int_{G\times H} |\varphi(t+s)  f(h + k') (1 - 1_{Q_K}(h) ) | \dd |\rho|(t,h)\leq
\varepsilon/3 \] 
for every $s\in G$ and $k'\in K$. On the other hand, by
continuity of $f$ and $\varphi$ and by compactness of $Q_K$, there obviously
exists an $n_1$ with
\[ \int_{G\times H} \left| \varphi (t+s_n)f(h + k_n)  - \varphi(t+s)  f(h+k) \right| 
\, 1_{Q_K}(h) \dd |\rho|(t,h)\leq
\varepsilon/3\] 
for all $n\geq n_1$. Putting this together, we infer
\[|\mu'(s_n,k_n) (\varphi) - \mu'(s,k)(\varphi)| \leq \varepsilon \]
for all $n\geq \max\{n_0,n_1\}$. 
\end{proof}

\begin{lemma} \label{factor} 
Let $\mcp$ be a measure cut and project scheme and $f : H\longrightarrow 
\mathbb C$ be admissible and continuous.
Then,  the dynamical system $(\varOmega(\nu_f),\alpha)$ is a factor of
  $(\TT,\beta)$ with factor map
\[ \mu   : \TT \longrightarrow \varOmega(\nu_f), \;\:\mu  ([s,k])
:=\mu'(s,k).\]
\end{lemma}
\begin{proof} 
  By Lemma \ref{property_mu} , $\mu' : G\times H \longrightarrow \MTB$ is
  continuous and $\widetilde{L}$-invariant. Thus, $\mu : \TT \longrightarrow
  \MTB$ is well defined and continuous. By definition, $\alpha_t (\mu ([s,k]))
  =\mu (\beta_t ([s,k]))$.  Thus, it remains to show that
\[\varOmega(\nu_f)=\mu   (\TT).\]
As $\mu  $ is a factor map, we have 
\begin{equation}\label{stern}
\mu  (\beta_t [0,0])  = \alpha_t
 (\mu  ([0,0]))=\alpha_t (\mu'(0,0)) =\alpha_t (\nu_f).
\end{equation}
By minimality of $\beta$ we have 
\[\TT = \overline{ \{\beta_t ([0,0]) : t\in G\} }.\]
Now, continuity of $\mu  $, compactness of $\TT$ and \eqref{stern}
imply
\[ \mu  ( \TT) = \mu  (\overline{ \{\beta_t ([0,0]) : t\in
  G\} })=\overline{ \{\alpha_t (\nu_f) : t\in G\} } =  \varOmega(\nu_f).\]
This finishes the proof.  \end{proof}

After these preparations we are ready to prove our first main result.

\begin{proof}[Proof of Theorem \ref{main}] 

(a) / (b) By Lemma \ref{factor}, $(\varOmega (\nu_f),\alpha)$ is a factor of
$(\TT,\beta)$ with factor map $\mu $.  Thus, by Theorem \ref{SAP}, $\nu_f$ is
strongly almost periodic and $(\varOmega (\nu_f),\alpha)$ carries the desired
group structure.

\smallskip

(c) As   $(\varOmega(\nu_f),\alpha)$ is a factor of $(\TT,\beta)$, it inherits
  spectral properties according to Fact \ref{thm_factor}. Now, $(\TT,\beta)$
  is well known to be uniquely ergodic and minimal with pure point spectrum
  and continuous eigenfunctions (see e.g. \cite{Martin2}). 
 As   $(\varOmega(\nu_f),\alpha)$ is a factor of $(\TT,\beta)$, the
  eigenvalues of $(\varOmega(\nu_f),\alpha)$ are eigenvalues of $(\TT,\beta)$
  as well. The eigenvalues of  $(\TT,\beta)$  can be determined easily
  \cite{Martin2}. Namely, each $\lambda\in \widehat{\TT}$ is an continuous
  eigenfunction to the eigenvalue $\lambda \circ \iota$, and this is a complete
  set of eigenfunctions. Now, the statement on the eigenvalues follows. 
\end{proof}

We finish this section with a proof of almost periodicity of $\nu_f$ for $f\in
C_c (H)$. While this statement is clear from the main theorem and the
abstract tools used above, it is instructive to give a direct proof.

\begin{lemma}\label{ap} Let $f\in C_c (H)$ be given. Then, $\nu_f$ is strongly
almost periodic. 
\end{lemma}
\begin{proof} Let $\varepsilon >0$ and $\varphi\in C_c (G)$ be arbitrary. We
  have to show that the set  $P$ of all $p\in G$   with
\[\| \delta_p \ast (\nu_f\ast \varphi) -  \nu_f \ast \varphi \|_\infty\leq \varepsilon\]
is relatively dense in $G$.  As $f\in C_c (H)$ and $\rho$ is
translation bounded, there exists an open neighbourhood $V$ of $0\in H$ with
\[\left|\int_{G\times H}   \varphi( t - s )   f(h + k)  \dd \rho(s,h) -
\int_{G\times H}  \varphi(t - s ) f(h) \dd \rho(s,h)\right|\leq \varepsilon\]
for all $t\in G$ and $ k\in V$. 
By $\widetilde{L}$-invariance of $\rho$, we have
\[ (\nu_f \ast \varphi) (t - \ell) = \int_{G\times H} \varphi (t -\ell -s) f(h) \dd \rho(s,h) =
\int_{G\times H}\varphi ( t - s)   f(h + \ell^\star)  \dd \rho(s,h)\]
for all $(\ell,\ell^\star)\in \widetilde{L}$. 
Putting the last two equations together, we see that 
$P$ contains all $\ell\in L$ with $\ell^\star
\in V$. As $V$ is open, this set is relatively dense. 
\end{proof}

\section{Proof of Theorem \ref{main_autocorrelation}}

In this section we provide a proof of Theorem
\ref{main_autocorrelation}. We need a preparatory result.

\begin{prop}\label{fubini} Let $\Oomega$ be a measure dynamical 
system with invariant measure $m$.  For arbitrary $\varphi,\psi 
\in C_c (G)$
$$ \int_\varOmega \int_G \int_G |\varphi (s+t) \psi (s)| \dd 
|\omega|(s) \dd |\widetilde{\omega}| (t) \dd m (\omega)< \infty.$$
\end{prop}
\begin{proof} This follows easily from uniform translation 
boundedness of $\omega\in \varOmega$. 
\end{proof}

\begin{proof}[Proof of Theorem \ref{main_autocorrelation}]  By Proposition
  \ref{integrability}, the function $f$ is integrable. By Theorem \ref{main}, 
$(\Omega(\nu_f)),\alpha)$ is uniquely ergodic. Denote the unique invariant
measure by $m$. Fix $\omega\in\Omega(\nu_f)$.

\smallskip

Recall now the definition of the measure $\gamma = \gamma_m$ as 
$$\gamma (\varphi) := \int_\varOmega \int_G \int_G \varphi ( s + t) \psi
(t) \dd \omega(s) \dd \widetilde{\omega}(t) \dd m (\omega)$$
for $\varphi \in C_c(G)$, where $\psi\in C_c (G)$ is arbitrary with 
$\int \psi(t) \dd t =1$.

\smallskip

Define $F$ on $\varOmega$ by $F(\omega) := \int_G \int_G \varphi (s+t)
\psi (t) \dd\omega(s) \dd\widetilde{\omega} (t).$ By Proposition
\ref{thm_factor} (a), we then have
$$\gamma(\varphi) = \int_\TT F(\mu(\xi)) \dd\xi.$$ 
Let $Z$ be a fundamental cell of $\widetilde{L}$ in $G\times H$.  In order to avoid a tedious factor  $1 / (m_G \times m_H) (Z)$ in the subsequent discussion, we will assume without
loss of generality that $(m_G \times m_H) (Z) =1$.

\smallskip

Recalling  $\mu'(s,h) = \mu([s,h])$ and applying the discussion before
Theorem 2, to $\sigma = m_G\times m_H$ instead of $\rho$, we obtain
$$\gamma (\varphi) = \int_Z F(\mu' (r,v)) \dd (m_G \times m_H) (r,v).$$
Unwinding the definitions then gives
$$\gamma (\varphi) = \int_Z \int_{G\times H} \psi (- t -r)
\left(\int_{G\times H} \varphi (s - t) f(h + v) \dd \rho (s,h)\right)
\overline{f} (k+v) \dd\overline{\rho} (t,k) \dd(m_G \times m_H) (r,v).$$
Note that the argument of $\varphi$ does not include
an $r$, as integration over $\omega$ contributes an $r$ and
integration over $\widetilde{\omega}$ contributes an $-r$ to the
argument of $\varphi$. We now use the $\widetilde{L}$-invariance of $\overline{\rho}$ to obtain
$$ \gamma(\varphi) = \int_Z \sum_{(l,l^\star)\in \widetilde{L}}
\left(\int_Z \psi (-t -l -r) G(t + l,v) \overline{f} (k + v + l^\star)
\dd\overline{\rho} (t,k) \right) \dd(m_G \times m_H) (r,v)$$ with $G(t,v)
= \int_{G\times H} \varphi (s -t) f(h +v) \dd\rho (s,h)$. By
$\widetilde{L}$-invariance of $\rho$, we have $G(t +l, v) =
G(t,v+l^\star)$, and we infer
$$\gamma(\varphi) = \int_{G\times H} \int_Z \psi (-t - r) G(t,v)
\overline{f} (k + v) \dd\overline{\rho} (t,k) \dd (m_G \times m_H)
(r,v).$$ By Proposition \ref{fubini}, we can now interchange the order
of integration. Carrying out the integration over $m_G (r)$ gives $1 =
\int_G \psi (-t -r ) \dd m_G (r)$, by assumption on $\psi$. The
integration over $m_H (v)$ yields $(f \ast \widetilde{f}) (h -k)$.
Altogether, we end up with
$$\gamma(\varphi) = \int_{G\times H} \left (\int_Z  (f \ast \widetilde{f}) 
(h- k) \varphi (s-t) \dd\rho (s,h)\right) \dd \overline{\rho} (t,k).$$

We now split the integration over $G\times H$ into integrations on
translates of the fundamental cell, yielding
$$\gamma (\varphi) = \int_Z \int_Z \left(\sum_{(l,l^\star)\in
\widetilde{L}} (f\ast \widetilde{f}) (h -k - l^\star) \varphi (s-t -l)
\right) \dd\rho(s,h) \dd\overline{\rho} (t,k).$$ This is the desired
formula.

\smallskip

If $\rho = \delta_{\widetilde{L}}$, then $\rho_\TT =
\delta_{[0,0]}$, and $f$ is square integrable by Proposition
\ref{squareintegrability}. This yields the remaining statements.
\end{proof}

\section{Weyl theorem on uniform distribution}\label{Weyl}
In this section we provide a proof of a Weyl result on uniform
distribution, Theorem \ref{theorem_Weyl} and derive two corollaries.
Such a result is of independent interest (and well known for usual
model sets).  Moreover, one of the corollaries will be used later when
we calculate the Fourier Bohr coefficients by a limiting procedure.

Our proof is inspired by \cite{Hof2,Moody2001}, with one
modification: we realize that the functional $\Lambda$ defined below
is translation invariant and hence must be a multiple of Haar
measure. This allows us to get a grip on our rather abstract situation
without more effort than the mentioned works. 

\smallskip

Various boundary and non boundary type terms have to be estimated.  To
do so we use the following proposition (see \cite{LM} for similar
results as well).

\begin{prop} Let $V\subset  G$ be a nonempty, open relatively compact set. 
Let $C\geq 0$ and  $\nu \in \MCV$ be arbitrary. Then, 
$$ |\nu| (B) \leq \frac{m_G (B - V)}{m_G(V)} \, C$$
for every relatively compact $B\subset G$. 
\end{prop}

\begin{proof} A direct calculation shows
$\chi_{B} \leq \frac{1}{m_G (V)} \chi_{B - V} \ast \chi_V.$
This gives
\begin{eqnarray*}
|\nu| (B) &\leq &    \frac{1}{m_G (V)} \int \int \chi_V (t -s) \, 
\chi_{B - V} (s) \,{\rm d}s\,{\rm d}|\nu|(t)\\
(Fubini) \:\;&=&  \frac{1}{m_G (V)}\int \chi_{B-V} (s) \left(\int  
\chi_V (t -s){\rm d}|\nu| (t)\right) {\rm d}s \\
&\leq & \frac{m_G (B - V)}{m_G(V)} C.
\end{eqnarray*}
This finishes the proof.  \end{proof}

With this proposition we can easily derive the following ``uniform''
version of Lemma 1.1 of \cite{Martin2}. More precisely, Lemma 1.1
deals with a single translation bounded measure. Here, we consider all
of $\MCV$ simultaneously.

\begin{lemma}\label{hilfe}  
Let $V\subset G$ open, nonempty and relatively compact be given. 
Let $(B_n)$ be a van Hove sequence in $G$ and set
$$ D:=\sup_{n\in\mathbb N}\frac{m_G (B_n - V)}{m_G (B_n)} < \infty.$$

(a) For every $C\geq 0$, 
$$\sup_{\nu \in \MCV, n \in \NN} \frac{|\nu| (B_n)}{m_G (B_n)} < C \, \frac{D}{m_G (V)}. $$

(b) For every $C\geq 0$ and $K\subset G$ compact
$$\lim_{n\to \infty} \sup_{\nu\in \MCV} \frac{|\nu| (\partial^K B_n)}{m_G (B_n)} =0.$$

(c) For   every $C\geq 0$ and  every $\varphi\in C_c (G)$ 
$$ \lim_{n\to \infty} \sup_{\nu \in \MCV} \frac{1}{m_G(B_n)} \left|
\int_G (\varphi \ast \nu) (t)\, \chi_{B_n}(t) \,{\rm d}t - \left(\int_G \varphi (t)
\,{\rm d}t\right) \, \nu (B_n)\right| = 0.$$
\end{lemma}
\begin{proof} 

(a) This is immediate from the previous proposition. 

(b) Set $\widetilde{K}:=(K\cup\{0\}) - \overline{V}.$ Then,
$\widetilde{K}$ is compact and
$$ \partial^K B - V \subset \partial^{ \widetilde{K} } B$$ 
for every $B\subset G$. Here, the $K$- boundary of a set was
defined in \eqref{rand}.  Thus, the previous proposition gives
$$ |\nu|(\partial^K B ) \leq m_G ( \partial^{ \widetilde{K} } B)
\frac{C}{m_G (V)},$$ 
and the statement follows as $(B_n)$ is a van Hove sequence.

(c) A short calculation (see proof of Lemma 1.1 (c) in \cite{Martin2}
as well) gives
$$ \left| \int_G (\varphi \ast \nu) (t)\, \chi_{B_n}(t)\,{\rm d}t - 
\left(\int_G \varphi (t) \,{\rm d}t\right) \, \nu (B_n)\right| 
\leq \int_G |\varphi(t)| \,{\rm d}t \,
|\nu| (\partial^K B_n),$$
and the statement follows from (b).

\end{proof}

\begin{theorem}[Weyl Theorem]\label{theorem_Weyl} 
Let  a measure cut and project scheme $\mcp$ be
  given. Let $f\in C_c (H)$ be arbitrary. Then, 
$$\lim_{n\to \infty} \frac{1}{m_G(B_n)} \mu(\xi) (\chi_{B_n}) =
\frac{\rho_\TT (1)}{(m_G\times m_H)_\TT (1)} \int_H f(h) \,{\rm d}h$$
uniformly in $\xi \in \TT$, where $\rho_\TT$ is defined in equation 
\eqref{form:rt}. 
\end{theorem}

\begin{proof} 
Note that the map $\mu : \TT \longrightarrow \MTB$  depends on 
$f$.  As shown in Lemma \ref{property_mu}, it is continuous. Thus, 
$\mu_\varphi : \TT \longrightarrow \CC$, $\xi\mapsto\mu(\xi)(\varphi)$, 
is continuous for every $\varphi \in C_c (G)$. Choose $\varphi 
\in C_c (G)$ with $\int_G \varphi (s) \dd s =1$. 
Define $\check{\varphi}$ by $\check{\varphi} (t) :=\varphi (-t)$ and note 
$\int \check{\varphi} \dd s =1$. 

As $(\TT,\beta)$ is uniquely ergodic, and $(B_n)$ is a van Hove sequence, 
the limit
$$\lim_{n\to \infty} \frac{1}{m_G(B_n)} \int_G \chi_{B_n} (t) 
\mu_\varphi(\beta_t \xi) \, {\rm d}t$$
exists uniformly in $\xi  \in \TT$.
A short calculation shows  $\mu_\varphi 
(\beta_t \xi) = (\check{\varphi} \ast \mu(\xi)) (t)$, and we infer from part 
(c) of the previous lemma that the limit
$$\Lambda (f) := \lim_{n\to \infty} \frac{1}{m_G(B_n)} \mu(\xi)
(\chi_{B_n})$$
exists uniformly in $\xi\in \TT$.  Apparently, $\Lambda : C_c (H) 
\longrightarrow \CC$ is a linear functional. 

\smallskip

We next show that $\Lambda$  has a certain boundedness property: 
Choose $K\subset H$ compact. Let $V\subset G$ be open, nonempty and 
relatively compact. As $\rho$ is translation bounded (see Proposition 
\ref{rho_tb}), 
$$ A:=\sup\left\{\int \chi_{V} (s+t) \chi_K (h +k) \,{\rm d}|\rho|(s,h) : 
(t,k)\in G\times H\right\} < \infty.$$
Thus, $\varOmega(\nu_f)\subset \cM_{A \|f\|_\infty,V}$ whenever 
$\supp f\subset K$. Therefore, part (a) of the previous lemma gives 
the existence of a constant $C_K$ such that 
\begin{equation}\label{star}
|\Lambda(f)| \leq C_K \|f\|_\infty
\end{equation}
for every $f\in C_c (H)$ with $\supp f \subset K$. This shows that $\Lambda$ 
is a measure on $H$. 

\smallskip

It is not hard to see that $\Lambda (f(\cdot - l^\star)= \Lambda(f)$ for
every $f\in L^\star$. As $L^\star$ is dense in $H$, we infer from the
continuity property \eqref{star} that $\Lambda (f(\cdot -h)) =
\Lambda(f)$ for every $h\in H$.

\smallskip

To summarise, $\Lambda$ is a translation
invariant measure on $C_c (H)$. Thus, $\Lambda$ is a
multiple of the Haar measure and  there exists $c_\rho\in \CC$ with
$$\Lambda(f) = c_\rho \int_H f(h) \,{\rm d}h.$$

In order to determine $c_\rho$, we choose a van Hove sequence $(C_m)$
in $H$. Let $\psi \in C_c (H)$ with $\int_H \psi \dd h=1$ be given and
set $f_m :=\chi_{C_m}\ast \psi$. Then, $f_m$ is a smoothed version of
$\chi_{C_m}$ and, in particular,
$$ c_\rho =  \lim_{m\to \infty} \frac{1}{m_H (C_m)} \Lambda (f_m). $$
Consider for fixed (and large) $m\in \NN$
$$\frac{1}{m_H (C_m)} \Lambda (f_m)=\lim_{n\to \infty} \frac{1}{ m_H(C_m) 
m_G (B_n)} \int_{G\times H} f_m (h) \chi_{B_n} (s) \,{\rm d}\rho (s,h).$$
For $n$ and $m$ large $  \int_{G\times H} f_m (h) \chi_{B_n} (s) \,{\rm d}\rho
(s,h)$  is, up to a boundary term, equal to
$$ \rho_\TT (1) \, \sharp \left\{ \widetilde{L}\cap (B_n \times C_m
  )\right\},$$ 
and $m_H(C_m) m_G (B_n) = (m_G \times m_H) (B_n \times C_m) $ is, up to a
boundary term, equal to
$$ (m_G \times m_H)_\TT (1) \, \sharp \left\{\widetilde{L}\cap (B_n \times
C_m)\right\},$$ where $\sharp$ denotes the cardinality. This easily gives the
desired value of $c_\rho$.
\end{proof}

\begin{coro} \label{Weyl_admissible} 
Let a measure cut and project scheme $\mcp$ be given. Let $f : H
\longrightarrow \CC$ be continuous and admissible. Then,
$$\lim_{n\to \infty} \frac{1}{m_G(B_n)} \mu(\xi) (\chi_{B_n}) =
\frac{\rho_\TT (1)}{(m_G\times m_H)_\TT (1)} \int_H f(h) \,{\rm d}h$$
uniformly in $\xi \in \TT$. 
\end{coro}

\begin{proof} 
This follows by approximation. Let $(B_n)$ be a
van Hove sequence in $G$.  Let $V\subset $ be open, nonempty and
relatively compact, and choose $\varphi \in C_c (G)$ with $0\leq \chi_{V}
\leq \varphi$. For every $\varepsilon >0$ we can then find by
admissibility a $Q'\subset H$ compact such that
$$ \int_{G\times H} |f(h +k ) (1 - \chi_{Q'} (h +k)) \varphi( s
+t)|\,{\rm d}|\rho|(s,h)\leq \varepsilon$$ for all $(t,k)\in G\times
H$.  Thus, the measures $\varphi \mapsto \int_{G\times H} f(h) (1 -
\chi_{Q'} (h)) \varphi( s)\,{\rm d}\rho(s,h)$ belong to $\MCV$ with
arbitrarily small $C$, provided $Q'$ is chosen large enough.  By  the explicit dependence on $C$ in (a) of
Lemma \ref{hilfe}, we can then find $Q_1\subset H$ compact with
\begin{equation} \label{term1}
\frac{1}{m_G (B_n)} \int_{G\times H} |f(h)( 1 - \chi_{Q_1} (h))|
\chi_{B_n} (t) \,{\rm d} |\rho|(t,h) \leq \varepsilon
\end{equation}
for all $n\in \NN$. By Proposition \ref{integrability}, we can find $Q_2
\subset H$ compact
\begin{equation}\label{term3}
\frac{\rho_\TT (1)}{(m_G\times m_H)_\TT (1)} \int_{H} |f(1 -
\chi_{Q_2})|\,{\rm d}h \leq \varepsilon.
\end{equation}
 Set $Q:= Q_1\cup Q_2$ and choose $\chi\in C_c (H)$ with $\chi_{Q}\leq
\chi$.  We can now write $f = f \chi + f (1 - \chi)$ and set $T_\xi^n
(f) := \frac{1}{m_G(B_n)} \mu(\xi) (\chi_{B_n})$ and $ T(f)
:=\frac{\rho_\TT (1)}{(m_G\times m_H)_\TT (1)} \int_H f(h) \,{\rm d}h$.
Then 
$$ |T^n_\xi (f) - T(f)|\leq |T^n_\xi (f) - T^n_\xi (f \chi)| + |
T^n_\xi (f \chi) - T(f\chi)| + |T (f\chi) - T(f)|.$$ Now, the first
term is smaller than $\varepsilon$ by \eqref{term1} for all $n$, the
second term goes to zero for $n\longrightarrow \infty$ by Weyl's
theorem. The last term is smaller than $\varepsilon$ by \eqref{term3}.
As $\varepsilon>0$ is arbitrary, this gives the desired convergence
statement.
\end{proof}

\begin{coro}\label{Weyl_riemann}
Let a measure cut and project scheme $\mcp$ be given with $\rho \geq 0$. 
Let $f : H \longrightarrow \RR$  be  Riemann integrable. Then, 
$$\lim_{n\to \infty}  \frac{1}{m_G(B_n)} \mu(\xi) (\chi_{B_n}) =
\frac{\rho_\TT (1)}{(m_G\times m_H)_\TT (1)} \int_H f(h) \,{\rm d}h$$
uniformly in $\xi \in \TT$. 
\end{coro}
\begin{proof} Again, this follows by approximation. As $f$ is Riemann 
integrable, there exists for every $\varepsilon >0$ $\varphi,\psi \in C_c (H)$ 
with
$$ \varphi \leq f  \leq \psi\;\:\;\mbox{and} \int (\psi - \varphi) \,{\rm d} h \leq 
\varepsilon.$$
Now, by Weyl's theorem Corollary \ref{Weyl_admissible}, the desired convergence holds for both $\varphi$
 and $\psi$ and the corollary follows easily. 
\end{proof}

\section{Fourier-Bohr coefficients and the proof of Theorem
  \ref{main_diffraction}} \label{Fourier}

In this section we provide a proof of Theorem
\ref{main_diffraction}. Throughout, we will assume that a measure cut
and project $\mcp$ and an admissible continuous $f:H\to\mathbb C$ is given. 

\begin{lemma}\label{gute_lambda} 
Let $\lambda\in \widehat{\TT}$ and $\xi \in \TT$ be given. Then, the limit
$$
c_\lambda (\xi) := \lim_{n\to \infty}  \frac{1}{m_G(B_n)} \mu (\xi) (\chi_{B_n}
\overline{\lambda \circ \iota})
$$
exists and 
$$ 
c_\lambda (\xi) = \lambda(\xi) \frac{\rho_\TT (\lambda)}{(m_G\times
  m_H)_\TT (1)}  \int f(h) \, (\lambda \circ \kappa)(h) \,\dd h.
$$
\end{lemma}

\begin{proof} 
This follows from a direct calculation using Weyl's Theorem. Choose  
$\xi = [s,k]\in  \TT$ arbitrary. Using $\lambda(\xi) = (\lambda\circ 
\iota)(s) \cdot (\lambda \circ \kappa) (k)$, we obtain after a short 
calculation
$$ 
\mu(\xi)(\chi_{B_n} \overline{\lambda \circ \iota}) = \lambda(\xi) 
\int\chi_{B_n} (s +t) (\lambda\circ \kappa) (h + k) f (h + k)
\overline{\lambda([t,h])} \,\dd \rho (t,h).
$$ 
We can now appeal to Weyl's theorem with
$f \, \lambda\circ \kappa$ instead of $f$ and $\lambda([\cdot
])\, \rho$ instead of $\rho$.
\end{proof}

\begin{lemma}\label{punktanteil} 
Let $\nu$ be a translation bounded measure on $G$ and $\gamma$ 
the associated autocorrelation. Let $\sigma\in \Ghat$ be given. If
$$
c_\sigma= \lim_{n\to \infty}  \frac{1}{m_G(B_n)} \nu 
(\chi_{B_n}\overline{\sigma})
$$
exists for every van Hove sequence $(B_n)$, then
$$
\gammahat(\{\sigma\}) = |c_\sigma|^2.
$$
\end{lemma}
\begin{proof} 
For $G=\RR^d$, this is proven by Hof in \cite{Hof}. The proof 
can be adapted to our more general situation.  More precisely, 
as discussed in section \ref{difftheo}, we have 
$\gamma = \lim_{n\to\infty} \frac{1}{m_G(B_n)} \nu_{B_n} \ast 
\widetilde{\nu_{B_n}}$. Thus,  Lemma \ref{hilfe} shows that 
there exists $C>0$ and $V\subset G$ open and relatively compact, 
with
\begin{equation}\label{assumption2}
\frac{1}{m_G(B_n)} \nu_{B_n} \ast \widetilde{\nu_{B_n}}\in
\MCV,\,\mbox{for all $n\in\NN$ and} \, \gamma\in\MCV.
\end{equation}
Moreover, by Theorem
11.3 of \cite{GdeL}, we have
\begin{equation}\label{assumption1}
\gammahat(\{\sigma\}) = \lim_{n\to \infty} \frac{1}{m_G (B_n)} \gamma
(\chi_{B_n} \overline{\sigma}).
\end{equation}
Given \eqref{assumption1} and \eqref{assumption2}, we can conclude 
as in the proof of Theorem 3.4 of \cite{Hof}.
\end{proof}



We can now come to the proof of Theorem \ref{main_diffraction}.

\begin{proof}[Proof of Theorem \ref{main_diffraction}]
By Theorem \ref{main}, $(\varOmega(\nu),\alpha)$ has pure point
dynamical spectrum with eigenvalues contained in $\{\lambda\circ 
\iota:\lambda\in \widehat{\TT}\}$. Thus, Theorem \ref{aequivalenz} 
gives that $\gammahat$ is a pure point measure which can be 
written in the form
$$
\gammahat = \sum_{\lambda\in \widehat{\TT}} w_{\lambda}
\delta_{\lambda\circ\iota}
$$
with suitable $w_\lambda$, $\lambda\in\widehat{\TT}$. For these 
$w_\lambda$, we obtain from Lemma \ref{gute_lambda} and Lemma 
\ref{punktanteil}
$$ 
w_\lambda = |c_{\lambda}(\xi)|^2,
$$
where the right hand side does not depend on $\xi\in\TT$. 
\end{proof}

\noindent \textbf{Remark.} In the situation discussed in this 
section, it is possible to show that
$$
\lim_{n\to \infty}  \frac{1}{m_G(B_n)} \mu (\xi) (\chi_{B_n} 
\hat{s}) =0
$$
for all $\xi \in\TT$, whenever $\hat{s}\neq \lambda \circ \iota$ 
for a $\lambda\in \widehat{\TT}$, see \cite{LM}. 

\section{Dense Dirac combs}\label{Dense} 

In this section, we restrict to the situation of $G=\mathbb R^d$, 
$H=\mathbb R^m$ and $\rho=\delta_{\widetilde L}=\sum_{x\in\widetilde
L}\delta_x$. This setup (with $\widetilde L=\mathbb Z^{d+m}$) appears 
in \cite{Hof} by regularising characteristic functions of model set 
windows. In a more general framework, the above setup is analysed in 
\cite{Ric}. Both situations are subsumed by our theory, as we will 
now show. Note that \cite{Ric} assumes a cut-and-project scheme 
with the additional assumption that $L$ is dense in $G$ and that 
the projection $\pi_{\rm int}$ is one-to-one between $\widetilde L$ 
and $L^\star$. We will not need these assumptions. 

Our arguments rest on the following special case of Corollary 
\ref{Weyl_riemann}, known as the density formula (see \cite{Martin1} 
as well). In order to formulate it in the variant discussed in 
\cite{Martin1,Ric}, we introduce the notation
$$\oplam(W) :=\{x\in L : x^\star \in W\}$$
for $W\subset H$ relatively compact. Thus, 
the connection with the translation bounded measures discussed so far 
is given by the formula
$$ \nu_{\chi_W} =  \sum_{x\in \oplam(W)} \delta_x.$$

\begin{coro}[\bf Density formula] 
Let a cut and project scheme $(G,H,\widetilde L)$ be given. 
Let $W\subset H$ be relatively compact such that $\chi_W$ is 
Riemann integrable. Then, for every van Hove sequence $(B_n)$
in $G$,
\begin{equation*}
\lim_{n\to\infty} \frac{1}{m_G(B_n)}
\left( \sum_{x\in\oplam(W+u)\cap (B_n+s)}1\right)
=\frac{m_H(W)}{(m_G\times m_H)_{\mathbb T}(1)}
\end{equation*}
uniformly in $s\in G$ and in $u\in H$.
\qed
\end{coro}

\smallskip

According to Proposition \ref{char_admissible}, a continuous function
$f:\mathbb R^m\to\mathbb C$ is admissible if and only if for arbitrary 
$\epsilon>0$ and $\varphi\in C_c(\mathbb R^d)$ there exists a compact 
$Q\subset\mathbb R^m$ with
\begin{equation*}
\sum_{x\in L}\left| \varphi(x+s)f(x^\star)\right|\left(1-1_Q(x^\star)\right)
\le\epsilon
\end{equation*}
for all $s\in\mathbb R^d$, where $1_Q$ denotes the characteristic function of
$Q$.

\begin{theorem}
Assume a cut-and-project scheme $(\mathbb R^d,\mathbb R^m,\widetilde L)$. Let
the function $f:\mathbb R^m\to\mathbb C$ be continuous, with
$|x|^{m+\alpha}|f(x)|\le C$ for all $x\in\mathbb R^m$, for some constants 
$C>0$ and $\alpha>0$. Then $f$ is admissible.
\end{theorem}

\noindent {\bf Remark.} The functions $f$ considered in \cite{Ric} 
satisfy the more restrictive condition $|x|^{m+1+\alpha}|f(x)|\le C$ for 
some constants $C>0$ and $\alpha>0$.

\begin{proof}
Let $(B_n)$ be a van Hove sequence in $\mathbb R^d$. For $l\in\mathbb N$, 
let $Q_l\subset\mathbb R^m$ denote the compact cube of sidelength $l$ 
centred at the origin. The density formula yields for $n>n_0$ large enough
the estimate
\begin{displaymath}
\left(\sum_{x\in\oplam(Q_1+u)\cap(B_n+s)}1\right)
\le 2\frac{m_G(B_n)}{(m_G\times m_H)_{\mathbb T}(1)}
\end{displaymath}
uniformly in $s\in\mathbb R^d$ and in $u\in\mathbb R^m$, since 
$m_H(Q_1)=1$. As $Q_{2(l+1)}\setminus Q_{2l}$ may be built from 
$(2l+2)^m-(2l)^m$ translated copies of $Q_1$, we obtain for $l\in\mathbb N$
\begin{displaymath}
\left(\sum_{{x\in(B_n+s)\cap L \atop x^\star \in Q_{2(l+1)}\setminus Q_{2l}}} 1\right)
\le 2\frac{m_G(B_n)}{(m_G\times m_H)_{\mathbb T}(1)} 
\left((2l+2)^m-(2l)^m\right)\le 2^{2m+1}\frac{m_G(B_n)}{(m_G\times m_H)_{\mathbb T}(1)}l^{m-1} 
\end{displaymath}
uniformly in $s\in\mathbb R^d$. The first estimate uses uniformity in $u$. 
To check admissibility of $f$, let $\epsilon>0$ and 
$\varphi\in C_c(\mathbb R^d)$ be given. Fix $n>n_0$ such that $\mbox{\rm supp}
(\varphi)\subset B_n$. We have for $l\in\mathbb N$ the estimate
\begin{displaymath}
\begin{split}
\sum_{x\in L}&\left| \varphi(x+s)f(x^\star)\right|\left(1-1_{Q_{2l}}(x^\star)\right) 
\le ||\varphi||_\infty \sum_{x\in(B_n-s)\cap L \atop x^\star \not\in Q_{2l}} \left| f(x^\star)\right|
= ||\varphi||_\infty \sum_{k=l}^\infty \sum_{x\in(B_n-s)\cap L \atop x^\star \in
  Q_{2(k+1)}\setminus Q_{2k}} \left| f(x^\star)\right|\\
& \le ||\varphi||_\infty \sum_{k=l}^\infty \frac{C}{k^{m+\alpha}}2^{2m+1}
 \frac{m_G(B_n)}{(m_G\times m_H)_{\mathbb T}(1)}k^{m-1}
= \widetilde C \sum_{k=l}^\infty \frac{1}{k^{1+\alpha}},
\end{split}
\end{displaymath}
where absolute convergence of the sum is used for reordering. The last
expression is a bound independent of $s\in\mathbb R^d$. Now choose 
$l\in\mathbb N$ such that the bound does not exceed $\epsilon$ and set 
$Q:=Q_{2l}$. By Proposition \ref{char_admissible}, we have shown 
that $f$ is admissible.
\end{proof}
According to Theorem \ref{main_autocorrelation}, we obtain for the 
autocorrelation
\begin{displaymath}
\gamma = \frac{1}{(m_G\times m_H)_\TT(1)}\sum_{l\in L} \eta(l)\,\delta_l, 
\qquad \eta(l)=\int_{\mathbb R^m} f(h)\overline{f(h-l^\star}) \dd h, 
\end{displaymath}
compare \cite[Thm.~9]{Ric}. For the diffraction formula, consider 
the {\it dual lattice} $(\widetilde{L})^\ast$ of $\widetilde L$, given by
\begin{displaymath}
(\widetilde{L})^\ast =\left\{(\gamma,\eta)\in 
\widehat{G}\times\widehat{H}\,:\,\gamma(l)\cdot\eta(l^\star)=1 
\mbox{ for all } (l,l^\star)\in\widetilde{L} \right\}.
\end{displaymath}
As $\widetilde{L}$ is a closed subgroup of $G\times H$, 
we infer that $(\widetilde{L})^\ast$ is a closed subgroup of $\widehat{G}
\times\widehat{H}$. By Pontryagin duality \cite[Thm.~4.39]{Fol}, the 
group $(\widetilde{L})^\ast$ is, as a topological group, isomorphic 
to $\widehat{\TT}$, an isomorphism from $\widehat{\TT}$ to 
$(\widetilde{L})^\ast$ being given by $\lambda\mapsto (\lambda 
\circ \iota,\lambda \circ \kappa)$ for $\lambda\in\widehat{\TT}$. 
Thus, $(\widetilde{L})^\ast$ is a discrete, cocompact subgroup of 
$\widehat{G}\times\widehat{H}$. Set now to $G\times H=\mathbb
 R^d\times\mathbb R^m$ with $\rho=\delta_{\widetilde L}$, and use the
canonical identification $\widehat{\mathbb R^n}\simeq\mathbb R^n$,
with $\eta(x)=e^{2\pi i\, \eta \cdot x}$ for $\eta\in\widehat{\mathbb R^n}$ and 
$x\in\mathbb R^n$. Then, Theorem \ref{main_diffraction} specialises to
\cite[Thm.~10]{Ric}.

\section{Injectivity of the factor map}\label{Injectivity} 

As discussed so far, a measure cut and project scheme $\mcp$ together
with an admissible continuous function $f$ yields a factor map $\mu :
\TT \longrightarrow \varOmega (\nu_f)$. In this section, we discuss
conditions to ensure that this factor map is one-to-one, i.e., an
isomorphism.

\medskip

First of all, note that injectivity of $\mu$ can be destroyed both 
by properties of $\rho$ and of $f$:

\medskip

\noindent \textbf{Example.} Let $\mcp$ be a measure cut and project scheme 
and $\rho :=m_G \times m_H$. Let $f$ be admissible and continuous.  Then,
$\mu(\xi) = (\int f \dd h)\,  m_G$ 
is independent of  $\xi \in \TT$, as can be seen by a direct calculation.

\medskip

\noindent \textbf{Example.} Let $\mcp$ be a measure cut and project 
scheme. Let $f$ be admissible and continuous with $f(\cdot - u) = f$ for some 
$u\neq 0$. Then, $\mu(\xi) = \mu(\xi + [0,u])$ for every $\xi\in 
\TT$, and $\mu$ is not injective. 

\medskip

Our aim in this section is to prove the following theorem. 

\begin{theorem}\label{thm_injectivity} Let a  cut and project scheme 
$(G,H, \widetilde{L})$ be given and  $f$ be admissible with respect to 
$\delta_{\widetilde{L}}$. Let $\rho$ be an $ \widetilde{L}$-invariant 
measure. If 
\begin{itemize}
\item $\rho_\TT (\lambda) \neq 0$ for every $\lambda \in \widehat{\TT}$, 
\item $f$ does not have a nontrivial period, i.e., $f(\cdot - u) \neq f$ 
for every $u\neq 0$,
\end{itemize}
then the map $\mu$ associated to $\mcp$  is one-to-one (and thus an 
isomorphism).
\end{theorem}

\noindent\textbf{Remark.}  (a)  Recall that admissibility of $f$ 
with respect to $\delta_{\widetilde{L}}$ implies admissibility of 
$f$ with respect to $\rho$, as shown in Lemma 
\ref{extension-of-admissibility}. \\
(b) Note that the first condition is
satisfied for $\rho = \delta_{\widetilde{L}}$. In this case $\rho_\TT
(\lambda) =1$ for every $\lambda \in \widehat{\TT}$.\\
(c) The condition that $f$ has no periods also appears in the context
of the usual cut and project schemes, where $f$ is a characteristic
function. There it is used to obtain a map from $\varOmega$ to $\TT$
(see \cite{Martin2,BLM}). In some sense our aim is similar. We use
this condition to prove injectivity of the map $\TT\longrightarrow
\varOmega$, which then also implies existence of a map from $\varOmega$
to $\TT$. 

\medskip

The proof will be given at the end of this section after a series 
of intermediate results.  We start by considering the case $\rho = 
\delta_{\widetilde{L}}$.

\begin{lemma}\label{spezial} Let a measure cut and project scheme 
$\mcp$ be given with $\rho = \delta_{\widetilde{L}}$.  If $f$ is 
continuous and admissible and does not have a nontrivial  period, 
then $\mu$ is injective. 
\end{lemma}

\begin{proof} As $\mu$ is a group homomorphism by Theorem \ref{main}, 
it suffices to show that $\mu(0) = \mu (\eta)$ implies $\eta =0$. 
Let $\eta=[(q,p)]$ be given with $\mu(0) = \mu (\eta)$. Thus,
\begin{equation}\label{hahaha}
\sum_{(l,l^\star)\in\widetilde L} f(l^\star) \delta_l = 
\sum_{(l,l^\star)\in\widetilde L} f(l^\star + p) \delta_{l + q}.
\end{equation}
As $L^\star$ is dense in $H$ and $f$ is continuous and does not vanish 
identically, there exists $l^\star\in L^\star$ with $f(l^\star)\neq 0$. 
Then, \eqref{hahaha} implies that there exists $l'\in L$ with 
$$ l = l' + q,$$
and we infer that $q = l - l' \in L$. Then, \eqref{hahaha} gives that 
$$ f(l^\star) = f (l^\star - q^\star + p)$$
for all $l\in L$. As $L^\star$ is dense in $H$ and $f$ is continuous, 
we obtain that $p - q^\star$ is a period of $f$. By assumption, all 
periods of $f$ are trivial, and we obtain $ p= q^\star$ and therefore
$$\eta =[(q,p)] = [(q,q^\star)] = 0.$$
This finishes the proof. 
\end{proof}

\begin{prop} 
Let a measure cut and project scheme $\mcp$ and an admissible
continuous $f: H\longrightarrow \CC$ be given. Let
$\mu:\TT\longrightarrow \varOmega(\nu_f)$ be the associated factor
map. Let $\varphi\in C_c (G)$ be arbitrary and define
$\mu_\varphi : \TT\longrightarrow  \CC, \,\xi\mapsto \mu(\xi)(\varphi).$
Then
$$ \widehat{\mu_\varphi} (\lambda) = \widehat{\varphi}(\lambda\circ
\iota)\,  \widehat{f} (\lambda\circ \kappa)\,  \widehat{\rho}(-\lambda).$$
\end{prop}
\begin{proof} 
Let $Z\subset G\times H$ be a fundamental cell of $\widetilde{L}$.
Then $G\times H = \widetilde{L} + Z$. For $\lambda\in \widehat{\TT}$ we can
then calculate $\widehat{\mu_\varphi} (\lambda)$ as follows:
\begin{eqnarray*} 
\widehat{\mu_\varphi}(\lambda) &=& \int_\TT \overline{(\lambda,\xi)} \mu_\varphi
(\xi) \,{\rm d}\xi\\ 
&=& \int_Z \overline{\lambda([s,k])} \left( \int_{G\times H} \varphi
(s + t) f (k+h) \,{\rm d}\rho(t,h)\right) \,{\rm d}s {\rm d}k\\ 
&=&\int_Z \overline{\lambda([s,k])}\left(\sum_{(l,l^\star)\in \widetilde{L}}  
\int_Z \varphi (s + l + t)f( k + l^\star + h) \,{\rm d}\rho(t,h)\right) 
\,{\rm d}s{\rm d}k \\ 
(Fubini)\;\:
&=& \int_Z \lambda([t,h]) \left( \sum_{(l,l^\star)\in \widetilde{L}} 
\int_Z \varphi (s + l + t) f( k + l^\star + h)\overline{\lambda([s,k]) \lambda([t,h])} \,{\rm d}s 
{\rm d}k\right) {\rm d}\rho(t,h)\\
&=&  \int_Z \lambda([t,h]) \left( \int_{G\times H} \varphi (s +t) f(k + h) 
\overline{\lambda([s +t, h +k])} {\rm d}s {\rm d}k\right) {\rm d}\rho (t,h)\\
&=& \int_Z \lambda([t,h]) \left( \int_{G\times H} \varphi (s) f(k)
\overline{ \lambda([s,  k])} \,{\rm d}s{\rm d}k\right) {\rm d}\rho (t,h)\\
&=& \int_Z \lambda([t,h])  \widehat{\varphi}(\lambda\circ
\iota) \widehat{f} (\lambda\circ \kappa)\,{\rm d}\rho (t,h)\\
&=& \widehat{\varphi}(\lambda\circ
\iota)\,  \widehat{f} (\lambda\circ \kappa)\, \widehat{\rho}(-\lambda).
\end{eqnarray*}
This finishes the proof. 
\end{proof}

We now come to the proof of injectivity.

\begin{proof}[Proof of Theorem \ref{thm_injectivity}] 
The proof of the general case will be reduced to the case treated in
Lemma \ref{spezial}. We want to show injectivity of the group
homomorphism $\mu$ associated to $f$ and $\mcp$. We will need as well
the group homomorphism $\mu^0$ associated to $f$ and
$(G,H,\widetilde{L},\delta_{\widetilde{L}})$. (Recall that $f$ is 
admissible with respect to $\delta_{\widetilde{L}}$ by assumption.)

If $\mu$ is not injective, there exists by Theorem \ref{main} an $\eta
\neq 0$ with $\mu (\xi + \eta) = \mu (\xi)$ for all $\xi \in
\TT$. Thus, $\mu_\varphi : \TT\longrightarrow \CC, \,\xi\mapsto
\mu(\xi)(\varphi)$ satisfies $\mu_\varphi (\xi + \eta) = \mu_\varphi
(\xi)$ for all $\xi\in \TT$ and $\varphi \in C_c (G)$. Taking Fourier 
transforms and using the previous proposition, we obtain
$$ 
(\lambda,\eta) \widehat{\varphi}(\lambda\circ \iota)\, \widehat{f}
(\lambda\circ \kappa)\, \widehat{\rho}(-\lambda) =
\widehat{\varphi}(\lambda\circ \iota)\, \widehat{f} (\lambda\circ
\kappa)\, \widehat{\rho}(-\lambda)
$$
for all $\lambda \in \widehat{\TT}$. As $ \widehat{\rho}(-\lambda) 
\neq 0$ for all $\lambda \in \widehat{\TT}$, this implies
$$
(\lambda,\eta) \widehat{\varphi}(\lambda\circ \iota)\, \widehat{f}
(\lambda\circ \kappa) = \widehat{\varphi}(\lambda\circ \iota)\, 
\widehat{f} (\lambda\circ\kappa)
$$
for all $\lambda \in \widehat{\TT}$. 
This, however, means that $(\eta,\lambda) \widehat{\mu^0_\varphi} 
(\lambda) = \widehat{\mu^0_\varphi} (\lambda)$ for all $\lambda \in 
\widehat{\TT}$, with $\mu^0_\varphi :\TT\longrightarrow \CC$ given 
by  $\mu^0_\varphi (\xi) = \mu^0(\xi)(\varphi)$. Taking the inverse 
Fourier transform, we obtain
$$
\mu^0_\varphi (\xi + \eta ) = \mu^0_\varphi (\xi)
$$
for all $\xi \in \TT$. As $\varphi \in C_c (G)$ is arbitrary, this 
gives 
$$
\mu^0 (\xi +\eta) = \mu^0 (\xi)
$$
for all $\xi \in \TT$. Now, $\mu^0$ is just the map treated in Lemma 
\ref{spezial}, where its injectivity is shown. Thus, we obtain 
$\eta =0$, and the proof is finished. 
\end{proof}

\section{A complementary result}\label{Complementary}

Let us shortly compare our results to the corresponding results for
the ``usual'' model sets. In our notation this case can be described
as follows (see Section \ref{Dense} as well): Let $\cp$ be a cut and 
project scheme and consider the measure 
\[ \rho := \delta_{\widetilde{L}}=\sum_{x\in \widetilde{L}} \delta_x\]
on $G\times H$. Now, let $\chi_W$ be the characteristic function of a 
compact set $W\subset H$, which is the closure of its interior and whose 
boundary has Haar measure $0$.  Then, we can form $\nu_{\chi_W}$ exactly 
as above. This measure has the form $\sum_{x\in \oplam(W)} \delta_x$ 
with the uniformly discrete set $\oplam(W) :=\{x\in L : x^\ast \in W\}
\subset G$. Identifying the measure $\nu_{\chi_W}$ with the uniformly 
discrete set $\oplam(W)$ in $G$, we can apply results of Schlottmann
\cite{Martin2} to obtain a factor map
\[\varPhi : \Omega(\nu_{\chi_W})\longrightarrow \TT\]
This map is $1:1$-almost everywhere, i.e. there exists a set $\TT_0$
of Haar measure zero in $\TT$ such that $\varPhi$ is one-to-one on
$\varPhi^{-1} (\TT \setminus \TT_0)$. It turns out that $\varPhi$ is
indeed not injective if $\nu_{\chi_W}$ is not periodic. More precisely, 
the following is proved in \cite{BLM}.

\bigskip

\noindent \textbf{Theorem.}\textit{
Let $\nu_{\chi_W}$ and $\varPhi$ be given as in the preceding paragraph. Then,
$\varPhi$ is injective if and only if $\nu_{\chi_W}$ is crystallographic, 
i.e. the set $\{t\in G : \alpha_t \nu_{\chi_W} = \nu_{\chi_W}\}$ is a cocompact
discrete subgroup of $G$.}

\bigskip

This theorem has the following consequence.

\begin{coro} 
Let the notation be as in the preceding theorem. If $\nu_{\chi_W}$ is not 
crystallographic, then $(\varOmega(\nu_{\chi_W}), \alpha)$ is 
not a factor of $(\TT,\beta)$. 
\end{coro}

\begin{proof} Assume the contrary. Then, there exists a factor map
\[\psi : \TT \longrightarrow \varOmega(\nu_{\chi_W}).\]
Thus, $\varPhi \circ \psi : \TT \longrightarrow \TT$ is a factor
map as well. With the canonical homomorphism $\iota : G\longrightarrow
\TT$, $\iota (t) = [t,0]$ and the definition of $\beta$ we therefore
obtain
\[ ( \varPhi \circ \psi ) (\iota(-t)) =( \varPhi \circ \psi ) ( \beta_t [0,0]) =  
\beta_t  ( \varPhi \circ \psi ) ([0,0]) = \iota(-t)\dotplus ( \varPhi \circ 
\psi ) ([0,0])=  ( \varPhi \circ \psi ) ([0,0]) \dotplus \iota(-t),\]
where we write $\dotplus$ to denote the product in $\TT$. As $\iota(G)$ is
dense in $\TT$ and $( \varPhi \circ \psi )$ is continuous, we infer
\[ ( \varPhi \circ \psi ) (\xi )= ( \varPhi \circ \psi ) ([0,0]) \dotplus \xi\]
for every $\xi \in \TT$.  In particular, $ ( \varPhi \circ \psi )$ is
injective as it is just translation by $( \varPhi \circ \psi ) ([0,0])$.  

On the other hand, by the previous Theorem, $\varPhi$ is
not injective as $\nu_{\chi_W}$ is not crystallographic.  As $\psi$ 
is a factor map, it is onto. Therefore, non-injectivity of $\varPhi$
leads to non-injectivity of $ \varPhi\circ \psi $. This
contradiction shows that $(\varOmega(\nu_{\chi_W}), \alpha)$ is not a factor
of $(\TT,\beta)$.
\end{proof} 

The corollary shows that our main result does indeed crucially depend
on the smoothness of the weight function $f$, as it becomes false for
characteristic functions of compact sets.  In some sense, the
dynamical systems associated with continuous weight functions are
closer to periodic systems than to model set systems.

\subsection*{Acknowledgements}  
The authors thank Michael Baake for important discussions and critical
comments on the manuscript.  A useful discussion with H.-J. Starkloff
on point measures is also gratefully acknowledged. D.L. would also
like to thank Peter Stollmann and Johannes Kellendonk for inspiring
conversations.  This work has been partially supported by the DFG.

\end{document}